\newcommand{\les}{\lesssim}

\def\pa{\partial}

\def\Db{{\bf {\bar {\,D}}}}
\def\Dh{{\bf {\hat {\,{D}}}}}
\def\Ds{{{\bf\cal D}}}

\def\laph{\hat{\Delta}}
\def\Null{\dot{\NN}^{-}}
\def\EE{{\cal E}}

\documentstyle[amssymb,amsfonts]{amsart}

\newenvironment{proof}{\noindent {\bf Proof} }{\endprf\par}
\def \endprf{\hfill  {\vrule height6pt width6pt depth0pt}\medskip}
\def\emph#1{{\it #1}}
\def\textbf#1{{\bf #1}}
\newcommand{\bea}{\begin{eqnarray}}
\newcommand{\eea}{\end{eqnarray}}
\def\beaa{\begin{eqnarray*}}
\def\eeaa{\end{eqnarray*}}

\def\ba{\begin{array}}
\def\ea{\end{array}}
\def\be#1{\begin{equation} \label{#1}}
\def \eeq{\end{equation}}

\newcommand{\nn}{\nonumber}

\def\Db{{^{(\LLb)}  D}}

\def\a{{\alpha}}
\def\b{{\beta}}
\def\ga{\gamma}
\def\Ga{\Gamma}
\def\de{\delta}
\def\De{\Delta}
\def\ep{\epsilon}
\def\eps{\epsilon}
\def\la{\lambda}
\def\La{\Lambda}
\def\si{\sigma}
\def\Si{\Sigma}
\def\om{\omega}
\def\Om{\Omega}

\def\th{\theta}
\def\ze{\zeta}
\def\nab{\nabla}

\def\Lb{{\underline{L}}}

\newcommand{\trchb}{\tr \chib}
\newcommand{\chih}{\hat{\chi}}
\newcommand{\chib}{\underline{\chi}}

\newcommand{\etab}{\underline{\eta}}
\newcommand{\chibh}{\underline{\hat{\chi}}\,}
\def\chih{\hat{\chi}}
\def\trch{\mbox{tr}\chi}
\def\tr{\mbox{tr}}

\def\NN{{\cal N}}

\def\FF{{\cal F}}

\def\JJ{{\cal J}}
\def\KK{{\cal K}}

\def\DD{{\cal D}}

\def\A{{\bf A}}
\def\B{{\bf B}}
\def\D{{\bf D}}
\def\F{{\bf F}}

\def\J{{\bf J}}
\def\M{{\bf M}}
\def\L{{\bf L}}

\def\R{{\bf R}}

\def\T{{\bf T}}
\def\g{{\bf g}}
\def\m{{\bf m}}

\def\LLb{{\bf \Lb}}

\def\RRR{{\Bbb R}}

\def\SSS{{\Bbb S}}

\def\B{{\bf B}}

\def\lap{\Delta}
\def\pr{\partial}

\def\c{\cdot}

\renewcommand{\div}{\mbox{div }}

\def\f14{{\frac{1}{4}}}
\def\f12{{\frac{1}{2}}}

\def\diag{{\mbox{diag}}}

\def\dual{{\,^\star \mkern-4mu}}
\def\2{{\overline 2}}

\begin{document}
\theoremstyle{plain}
  \newtheorem{theorem}[subsection]{Theorem}
  \newtheorem{conjecture}[subsection]{Conjecture}
  \newtheorem{proposition}[subsection]{Proposition}
  \newtheorem{lemma}[subsection]{Lemma}
  \newtheorem{corollary}[subsection]{Corollary}

\theoremstyle{remark}
  \newtheorem{remark}[subsection]{Remark}
  \newtheorem{remarks}[subsection]{Remarks}

\theoremstyle{definition}
  \newtheorem{definition}[subsection]{Definition}

\include{psfig}
\title[Kirchoff-Sobolev parametrix ]{A Kirchoff-Sobolev parametrix for  the
 wave equation and applications}
\author{Sergiu Klainerman}
\address{Department of Mathematics, Princeton University,
 Princeton NJ 08544}
\email{ seri@@math.princeton.edu}

\author{Igor Rodnianski}
\address{Department of Mathematics, Princeton University, 
Princeton NJ 08544}
\email{ irod@@math.princeton.edu}
\subjclass{35J10\newline\newline
The first author is partially supported by NSF grant 
DMS-0070696. The second author is partially 
supported by NSF grant DMS-01007791. Part of this work was done while he was visiting 
Department of Mathematics at MIT
}
\maketitle
\vspace{-0.3in}
\section{Introduction}
In this paper we propose a  construction of  
a   first order parametrix for solutions to the 
covariant,  tensorial  wave equation
\begin{equation}\label{eq:Wave}
\Box_\g \Psi = F,
\end{equation}
under minimum 
assumptions for the
 Lorentz manifold  $(\M,\g)$.
Here $\Psi$ and $F$ are $k$ tensor-fields 
on a     $3+1$ dimensional Lorentz manifold  $(\M,\g)$
and   
$$\Box_\g \Psi=\g^{\mu\nu} \D_{\mu}\D_{\nu}\Psi
$$ 
denotes
the covariant wave operator  on $\M$,  with $\D$ the Levi-Cevita 
connection defined by  $\g$. To  simplify the discussion below 
we  consider first  the scalar case 
\begin{equation}\label{eq:wave}
\Box_\g \psi = f.
\end{equation}
In Minkowski space $(\RRR^{3+1}, \m)$ 
with $\m=\diag\{-1,1,\ldots, 1\}$
 the wave operator on the left hand side of \eqref{eq:wave}
 is the standard
 D'Alembertian  $\Box=\m^{\a\b}\pr_\a\pr_\b$.
The general solution
of $\square \psi=f$ can be written in the form,
\be{eq:gen-sol-mink}\psi=\psi_f+\psi_0
\end{equation}
with $\psi_0$ a solution of the homogeneous equation
$\square \psi_0=0$ and  $\psi_f$ given by the 
  the Kirchoff formula,
\bea
\psi_f(t,x)&=&(4\pi)^{-1}\int_{0}^t\int_{|x-y|=t-s} |x-y|^{-1} f(s,y) ds
d\si(y)\nn\\ &=&(4\pi)^{-1}\int_{\RRR^{3+1}_+}\frac{1}{|x-y|}\de(t-s-|x-y|) f(s,y) ds
dy.\label{eq:Kirch1}
\eea
Here  $d\si(y)$ denotes  the area element of the sphere $|x-y|=t-s$ and
$\de$ represents the one dimensional Dirac measure supported at 
the origin. The homogeneous solution $\psi_0$ is fixed by initial data 
on the hyperplane $t=0$.

 One can also recast
 \eqref{eq:Kirch1} in the form
\be{eq:Kirch2}
\psi_f(t,x)=(4\pi)^{-1}\int_{\R^{3+1}_+} H(t-s)\de\big(-(t-s)^2+|x-y|^2\big)f(s,y) ds dy
\end{equation} 
where $H(t)$ is the Heavyside function
 supported on the positive  real axis
 and the expression  $|x-y|^2-(t-s)^2=d_0(p,q)^2$ is the square of
the  Minkowski distance function between  the vertex $p=(t,x)$ and the point
$q=(s,y)$ in
 the causal past $\JJ^{-}(p)\cap\R^{3+1}_+$ of the point  $p\in \R^{3+1}_{+}$. 
 All attempts
to extend  Kirchoff's formula  to a general  four
 dimensional curved space-time
are based on either  \eqref{eq:Kirch1} or \eqref{eq:Kirch2}. Thus the 
 first term in the  so called  Hadamard parametrix is constructed
 by replacing 
the  Minkowski distance function $d_0$  with the Lorentzian
distance function $d(p,q)$ defined by the metric $\g$.
 Thus one can set,
\be{eq:Kirch3}
\psi_f(p)=(4\pi)^{-1}\int_{\JJ^{-}(p)}r(p,q)\,\de\big( d^2(p,q)\big)\,f(q)\, dv(q)
\end{equation}
with $dv$ the volume element of the metric $\g$,  and $r(p,q)$
a correction factor which verifies a transport equation
 along  the null boundary of  $\JJ^{-}(p)$ and such that
  $r(p,p)=1$. The integral on the
right makes sense for the portion of  $\JJ^-(p)$ which belongs to a neighborhood
$\DD$  of $p$ where the geodesic distance function $d(p,q)$
is well defined and  sufficiently smooth. 
Typically one requires    $\DD$ to be
  causally geodesically convex, i.e. any two causally
separated  points  in    $\DD$ can be joined
by a unique geodesic in $\DD$. 
The local parametrix in $\DD$ is then defined 
 \be{eq:Kirch3'}
\psi_f(p)=(4\pi)^{-1}\int_{\JJ^{-}(p)\cap \DD}r(p,q)\,\de\big( d^2(p,q)\big)\,f(q)\, dv(q)
\end{equation}
The integral in \eqref{eq:Kirch3'} is supported on the portion of the 
 boundary $\NN^{-}(p)$ of $\JJ^{-}(p)$ included in $\DD$.

 The error term $\square_\g \psi_f-f$, however,   does not
vanish  unless $\g$ is the flat metric $\m$. 
One can improve   \eqref{eq:Kirch3}
by making successive corrections   based on solving a series of 
transport equations in $\JJ^{-}(p)\cap\DD$. In the process 
the error term can be made
 as smooth as we wish,  for given regularity of $f$,  at   the price
of requiring higher regularity of the metric $\g$, 
see \cite{Fried}. Moreover the resulting parametrix,
called Hadamard parametrix, is no longer supported just  on
the boundary of $\JJ^-(p)$.
One obtains a  solution of  \eqref{eq:wave} of the form,
\be{eq:hadamard}
\psi(p)=\int_{\JJ^{-}(p)\cap\DD} E_-(p,q) f(q)\, dv(q).
\end{equation}
with $E_-(p,q)=r(p,q)\de\big( d^2(p,q)\big)+\ldots$
is the retarded Green function of $\Box_\g$

The Hadamard parametrix \eqref{eq:hadamard}, which requires
 both  infinite smoothness of $\g$ and geodesic convexity
 for $\DD$ is ill suited 
for applications  to nonlinear problems. It turns
out that in many situations one does not need  
the precise representation \eqref{eq:hadamard} and  that 
in fact   the first order parametrix
 of type \eqref{eq:Kirch3}  suffices.   This fact  was  first
made use of  by  S. Sobolev,
see \cite{Sob},    to provide a proof of 
well-posedness for general
second order linear wave equations with variable
 coefficients. A similar
parametrix was later used by Y. C. Bruhat, see  \cite{Br}, in her famous local 
existence result for   the  Einstein vacuum equations. 
Both   \cite{Sob} and \cite{Br}  construct their first order parametrices,
which we  refer to as Kirchoff-Sobolev,  based on the  
flat space formula\footnote{It is easy to show that the two constructions
\eqref{eq:Kirch3} and \eqref{eq:Kirch4} differ in fact
 only by a normalization
factor at the vertex $p$.}
 \eqref{eq:Kirch1}. 
The generalization of \eqref{eq:Kirch1} to a curved space-time  proceeds 
from  the observation that  the function $u_p(s,y)= t-s-|x-y|$ is an optical
function, i.e.
\be{eq:opt1}
\m^{\a\b}\pr_\a u\pr_\b u=0,
\end{equation}
vanishing precisely on the past  null cone $\NN^{-}(p)$ 
with vertex at $p=(t,x)$ given by the equation $u_p=0$.
Letting $r=|x-y|$  one can easily check that 
\beaa
\square \big( r^{-1}\de(u_p)\big)&=&\big(\square \, r^{-1}\big)\,\de(u_p)
 +(-2 L(r^{-1})+r^{-1}\square u_p)\de'(u)\\
&+&\big(\m^{a\b}\pr_\a u_p\pr_\b u_p\big)\de''(u)=4\pi\de(p),\qquad 
\eeaa
with $\de(p)$ the four dimensional Dirac measure supported at $p$.
Indeed the terms  involving $\de''(u_p)$ and $\de'(u_p)$ both
 vanish, the first in view of \eqref{eq:opt1} and the second 
 because,
$$-2 L(r^{-1})+r^{-1}\,\square u_p=0,$$ with $L$ the null vectorfield along
$\NN^{-}(p)$ defined by $L=-\m^{\a\b}\pr_\b u_p \pr_\a$. On the other hand 
$\de(u_p)\square r^{-1}=\de(u_p)\lap r^{-1}=4\pi\de(p)$.

Based on this one can generalize
 \eqref{eq:Kirch1} to  a curved space-time by setting,
\be{eq:Kirch4'}
\psi_f(p)= \int_{\JJ^{-}(p)\cap\DD}  a(p,q)\,\de(u_p(q)) \,f(q)\, dv(q)
\end{equation}
where $u_p=u_p(q)$ is the  backward solution to the 
eikonal equation,
\be{eq:Eikonal}
\g^{\a\b}\pr_\a u\,\pr_\b u =0,
\end{equation}
vanishing on the past null cone $\NN^{-}(p)$,
and $a(q)=a(p,q)$ verifies the transport equation 
similar to that satisfied by $r^{-1}$ in flat space.
As in \eqref{eq:Kirch3}  we  need to restrict 
ourselves to a  neighborhood  $\DD$ of $p$  in which
solutions to \eqref{eq:Eikonal} remain smooth. 
 
 To explain  the restriction to the  neighborhood ${\cal D}$ to which the 
 integral in \eqref{eq:Kirch4'} is restricted we return for a moment to the 
 initial value problem in flat space-time. In the Minkowski space-time 
 model with the choice of an initial Cauchy hypersurface $\Si_0=\{t=0\}$
 the Kirchoff formula 
 \begin{equation}\label{eq:Kirch-flat}
 \psi_f(p)= (4\pi)^{-1} \int_{\JJ^-(p)\cap \JJ^+(\Si_0)}\frac 1{r(p,q)} {\de\left(u_p(q)\right)}
\, f(q)\,
 \,dv(q)
 \end{equation}
 with $p=(t,x), q=(s,y)$, $u_p(q)=t-s-|x-y|$ and $r(p,q)=|x-y|$, coincides at point $p$ with the solution 
 of $\Box \psi = f$ with zero initial data at $t=0$. The representation is valid for {\it any}
 point $p$ to the future of $\Si_0$ and the surface of integration 
 $$
 \NN^-(p)\cap \JJ^+(\Si_0)=\{ (s,y):\,\, t-s=|x-y|,\,\, s\ge 0\}
 $$
 is smooth with exception of the vertex point $p$. In a flat space-time model with the 
 Lorentzian manifold $\M=\R\times \Pi_a$, where $\Pi_a=\R^2\times \R/a {\bf Z}$ is a flat 
 cylinder of  ``width" $a$, the representation \eqref{eq:Kirch-flat} also coincides with 
 the solution of the inhomogeneous wave equation with zero initial data at $t=0$,
provided that we restrict ourselves to 
  points $p=(t,x)$ such that $t\le a$. For points $p=(t,x)$ with $t>a$  formula 
 \eqref{eq:Kirch4'} no longer\footnote{In fact the correct representation can be obtained
 by lifting the problem to the covering space $\R\times \R^2\times \R$, applying 
the Kirchoff formula
 and taking periodization in the last variable with the period $a$.} 
 represents  the solution of the inhomogeneous problem with 
 zero initial data at $t=0$. The null hypersurface $\NN^-(p)\cap \JJ^+(\Si_0)$ develops 
 singularities\footnote{Note that although past null geodesics intersecting, say at $t_*=t-a$ can be
 extended beyond $t_*$ they no longer belong to the boundary of the causal past of $p$.} 
 (scars) in the time interval $[0,t-a]$ due to intersecting null geodesics. This shows that the 
 accuracy of the Kirchoff formula in this case is restricted to the neighborhood 
 $\DD=\{(t,x):\,\,0\le t\le a\}$ of the Cauchy hypersurface $\Si_0$.

 To describe the situation in a general space-time $(\M,\g)$ we assume that 
 $\M$ is globally hyperbolic, i.e., there exists a Cauchy hypersurface $\Si\subset \M$ 
 with the property that each in-extendible past (future) directed causal curve from 
 a point $p$ to the future (past) of $\Si$ intersects $\Si$ once. We denote by 
 $\Si_+=\JJ^+(\Si)$ the future set of $\Si$. By finite speed of propagation the  solution 
$\psi(p)$
 of the wave equation $\Box_\g\psi=f$ at point $p\in \Si_+$ is completely determined by the 
 values of $f$ in $\JJ^-(p)\cap \Si_+$ and initial data for $\psi$ on $\JJ^-(p)\cap \Si$. 
 \begin{definition}
 We will say that $E_-(p,q)$ is the retarded parametrix for $\Box_\g$ at $p$ if
 $$
 \psi(p)=\int_{\JJ^-(p)\cap \Si_+} E_-(p,q) \, f(q)\, dv(q)
 $$
 coincides with the solution of the problem $\Box_\g\psi=f$ with zero initial data 
 on $\Si$. We will say that the first term in the expansion of 
 $E_-(p,q)$ -- distribution $\KK^-_p=a(p,q) \de \left(u_p(q)\right )$ -- is the retarded 
 Kirchoff-Sobolev parametrix.
 \end{definition}
 Let $\DD$ be a space-time  neighborhood of $\Si$. The expression  
\be{eq:Kirch4}
\psi_f(p)= \int_{\JJ^{-}(p)\cap\Si_+}  a(p,q)\de(u_p(q)) f(q) dv(q),\qquad p\in \DD
\end{equation}
is the Kirchoff-Sobolev approximation to the solution $\psi(p)$ of the wave equation 
$\Box_\g\psi=f$ with zero initial data on $\Si$. 
Clearly $\psi_f$ fails to be a solution to
\eqref{eq:wave} in the non-flat case. We write a  general solution
of \eqref{eq:wave} with zero initial data on $\Si$ in the form,
\be{eq:gen-wave}
\psi(p)=\psi_f(p)+{\cal E}_f(p)
\end{equation}
with ${\cal E}_f$ an error term.

In this paper we will:
\begin{enumerate}
\item Provide a careful  
derivation of \eqref{eq:Kirch4} and \eqref{eq:gen-wave}
for points $p$ in a  suitable  neighborhood $\DD$ of $\Si$
 and show that  the error term ${\cal E}_f$ can be 
expressed in the form,
\be{eq:error-f}
{\cal E}_f(p)=\int_{\JJ^{-}(p)\cap\Si_+}  {\cal E}(p,q)\,\de(u_p(q))\, \psi(q)\, dv(q)
\end{equation}
where  the smooth density ${\cal E}(p,q)$ depends only on geometric 
quantities associated to the null hypersurface $\NN^{-}(p)$.
\vskip 1pc
\noindent
We should note that classical constructions of the Kirchoff-Sobolev parametrix 
establish the error term ${\cal E}_f$ as explicitly dependent on the metric $\g$ and
its derivatives relative to some chosen system of coordinates. To our knowledge
the fact that  ${\cal E}_f$ is supported only on the boundary 
of the past set $\JJ^-(p)$ does not seem 
to have been fully  recognized and used in
applications. A similar observation was, prior
to this work,   communicated to us verbally  by V.~Moncrief. His claim, based on 
Friedlander's treatment of the Hadamard parametrix,  was the starting 
 point  of our  own investigations.
\vskip 1pc

\item Extend 
  formulas \eqref{eq:Kirch4} and \eqref{eq:error-f} to the covariant tensorial wave
equation \eqref{eq:Wave}.

\vskip 1pc
\noindent 
Once again the classical treatment of the tensorial wave equation introduces additional
coordinate dependent error terms. Our approach is entirely covariant.
\vskip 1pc
\item Provide a minimum set of conditions
for the local geometry of $\M$ near $p$ to ensure that the representation 
\eqref{eq:Kirch4} and \eqref{eq:gen-wave} holds true at $p$. We also 
make use of  our recent results from \cite{KR4} to show that for the Einstein vacuum 
space-times $(\M,\g)$,  with vanishing Ricci curvature, formulas 
\eqref{eq:Kirch4} and \eqref{eq:gen-wave} can be extended to points $p$ 
at  distance $t_*$ from $\Si$,  with $t_*$ dependent, essentially,  only on the $L^2$ norm 
of curvature\footnote{Note that classically the  construction of a Kirchoff-Sobolev
parametrix could only be justified   
for points $p$ such that $\JJ^-(p)\cap \Si$ belongs to a geodesically convex neighborhood
$\DD$ of $p$. As we note below this requires uniform
control for at lest two derivatives of the metric.} of $\g$. 

\item 
\vskip 1pc Our   formula  can be easily  adapted to gauge invariant  wave
equations. In section 4 of the paper we write down such a formula and show how it
can be used to give
 a very simple proof of the  Eardley-Moncrief  global existence
 result  for the Yang-Mills equation in the $3+1$ dimensional Minkowski space, 
see \cite{EM1},\cite{EM2}.
The remarkable fact about our approach is that it is entirely gauge independent;
we don't need to specify  any gauge condition\footnote{The method of 
\cite{EM1},\cite{EM2}
 was heavily dependent on the choice of a Cr\"onstrom gauge.    }.
\noindent 

\end{enumerate}

The size of the neighborhood $\DD$, mentioned above,  is first and foremost constrained by
the condition that the optical function $u$ is smooth. In the case of  a
Riemannian manifold the distance function from a point $p$ 
is smooth in a geodesically convex neighborhood of $p$ whose size
can be evaluated in terms of the $C^2$ norm of  the metric $\g$, as 
measured in a given system of 
coordinates. Alternatively, by a theorem of Cheeger, the size of this
 neighborhood depends
only on the pointwise bounds for the Riemann curvature tensor and a lower bound on 
the volume of a unit geodesic ball. For similar reasons 
the construction of  a  solution  $u_p$ to \eqref{eq:Eikonal}
is restricted to a  geodesically
convex neighborhood\footnote{Defined as the image of the exponential map 
$:T_p\M\to \M$ restricted to the largest convex subset 
of $T_p\M$ 
where it is a diffeomorphism.} 
$\DD$ of $p$.   Unlike the 
Riemannian case, however, a  purely geometric  characterization 
 of the size of a geodesically convex neighborhood of a point
 $p$ is not available and thus all known
parametrix constructions for wave equations
 had to be  restricted to  domains $\DD$ whose
size is
 determined by 
the $C^2$ norm of the metric $\g$ 
in a  given system of coordinates. Thus the Kirchoff-Sobolev representation 
would only hold for points $p$ at maximal distance $t_*$ from $\Si$ with 
$t_*$ dependent on the $C^2$ norm of the metric. As we shall explain below,
 such
demand on the regularity of the metric  would make the Kirchoff-Sobolev formula
impossible to apply to realistic nonlinear situations, such as Einstein's field
equations.

The importance of the classical  $C^2$ condition becomes apparent upon examining 
 the regularity of the   null boundary $\NN^-(p)$  
of the causal past $\JJ^{-}(p)$.
This set is   ruled by past  null geodesics $\ga(s)$ 
originating from $p$ and terminating at the points $\ga(s_*)$ beyond which one can 
find a time-like curve connecting $p$ and $\ga(s)$ with $s>s_*$, see \cite{HE}. Regularity
of $\NN^-(p)$  breaks down precisely at the terminal points $\ga(s_*)$. There are two
reasons for the  existence of a terminal point $\ga(s_*)$. 
\begin{enumerate}
\item $\ga(s_*)$ is a conjugate point.
\item $\ga(s_*)$ is a point of intersection of two different null geodesics.
\end{enumerate}
The existence of conjugate points is governed by the Jacobi equation for the Hessian  
$\D^2 u$ of  the optical function $u$,
$$
\D_\L (\D^2 u) + (\D^2u)^2 = \R(\,\c\,, \,\L,\,\c\,,\L)
$$
with $\L=-\g^{\a\b}\pr_\b u\pr_\a$ the null geodesic vectorfield 
along $\NN^{-}(p)$ and $\R$ the curvature tensor of $\g$.
 This formula indicates that, at least as far as the conjugate points are concerned,
 the terminal value of the affine parameter $s_*$ can be  bounded 
below by an upper bound on sectional curvature which, in turn,
 can be  controlled by a $C^2$  bound on the metric. 

 Uniform bounds of 
of the curvature tensor $\R$,
 or $C^2$ bounds for the metric $\g$, are however not very useful in applications to
nonlinear wave equations. For example in the  classical local  existence
result for the Einstein vacuum equations \cite{Br}, 
which is based on Kirchoff-Sobolev formula,  the $C^2$ requirement
is by itself worse\footnote{Additional losses of derivatives
lead to a $C^5$  result in \cite{Br}.}  when compared to
the  result in
\cite{HKM}  based on the Sobolev norm $H^s$, $s>5/2$.
It is for this reason alone that the Kirchoff-Sobolev 
parametrix has been abandoned in all rigorous work 
on nonlinear wave equations in favor of energy estimates
 and Sobolev inequalities.
The main goal of our paper is to  revive the 
Kirchoff-Sobolev parametrix by  constructing it and showing that 
in the particular case of the Einstein vacuum equations,
  $$\R_{\a\b}=0,$$
it is well-defined  under much  less stringent assumptions.
 For this task we rely in an essential way on  the results 
in \cite{Causal1}--\cite{Causal3} which show\footnote{Properly speaking
the results in \cite{Causal1}--\cite{Causal3} do not
 consider  the vertex $p$ yet the methods used
in those  papers can be  shown to extend 
to  cover the  case  of interest here. In fact, this forms the subject of the Q. Wang's 
thesis, Princeton University, 2006.} that the radius of
conjugacy along $\NN^{-}(p)$, expressed relative  to an
affine parameter of $\L$,  depends only on the size of 
the  geodesic 
flux of  curvature\footnote{This is an appropriate 
$L^2$ integral of the tangential components 
of the curvature tensor along $\NN^{-}(p)$, called curvature flux,
which will be defined below.}  $\FF_p$ along $\NN^{-}(p)$. These results are 
complemented by our recent  work \cite{KR4} where we 
establish the remaining part of a lower bound on the radius 
of injectivity  of $\NN^{-}(p)$, i.e., control of intersecting null geodesics from
 $p$, 
expressed relative to a  given time function. We achieve this 
by assuming, in addition to  the above mentioned bound on the curvature flux,  the
existence of a coordinate system
$x^\a$
 in $\DD$  relative to which the metric $\g$ is
 pointwise close to the flat Minkowski  metric. 

{\bf Acknowledgment.}\quad We would like to thank V. Moncrief for
fruitful discussions in connection with our work. He was
first to point out to us that a formula of type \eqref{eq:gen-wave}
with an error term  $\EE_f$ of the form \eqref{eq:error-f},  supported 
on the boundary of the past of $p$,  should hold 
true. His derivation,  based on  Hadamard's parametrix construction
 as formulated  in \cite{Fried},
differs  however significantly from ours.
We  would also like to point out that our invariant derivation 
of the Eardley-Moncrief  global regularity result for the $3+1$
dimensional Yang-Mills equations answers a
question first  raised  to us by him. 
\section{Basic definitions and main formula}
\subsection{Null cones}\label{se:null-cone}
Consider  a  spacelike hypersurface $\Si$,  a point $p$ to its future
 $\Si_{+}$ and $\JJ^{-}(p)$ its causal past. We start by  assuming 
the following local hyperbolicity  condition for the pair $(\Si,p)$:

{\bf A1.}  
{\it All past causal curves initiating  at points in a small  neighborhood  
of $\JJ^{-}(p)$ intersect $\Si$
at precisely one point. }

 Let  $\NN^{-}(p)$ be the  null boundary  of $ \JJ^{-}(p)$.  In general $\NN^{-}(p)$ is
an achronal, Lipschitz  hypersurface. 
 It is ruled by  the  null geodesics\footnote{
Every point in $\NN^{-}(p)\setminus\{p\}$ can be reached from $p$ 
 by a  past null geodesic  in $\NN^{-}(p)$.} from $p$,
corresponding to all past null directions in the 
tangent space $T_p\M$.  
These null geodesics 
can be parametrized by fixing  a future unit time-like vector
 $\T_p$ at $p$. Then,  for every direction $\omega\in \SSS^2$,
 with  $\SSS^2$ denoting the standard sphere in $\RRR^3$,  consider the  
 null vector $\ell_\omega$ in $T_p(\M)$,  
\be{eq:norm-geod}
\g(\ell_\omega, \T_p)=1,
\end{equation}
  and associate to it the  past null
geodesic $\ga_\omega(s)$ with initial data $\ga_\omega(0)=p$ and 
$\dot \ga_\omega(0)=\ell_\omega$. We can choose the parameter $s$ in 
such a way so that $\L=\dot\ga_\om(s)$ is geodesic.  Thus,
\be{eq:geod-null}
\D_\L\L=0,\qquad \g(\L,\L)=0,\qquad  \mbox{ and, at point  p,}\quad \g(\L, \T_p)=1
\end{equation}
As mentioned in the introduction the null cone  $\NN^{-}(p)$ 
is smooth as long as the  exponential map $(s,\om)\to \ga_\om(s)$
is a local diffeomorphism and no two geodesics,
corresponding to different direction  $\om\in \SSS^2$,
intersect. Thus for each $\om\in \SSS^2$ either $\ga_\om (s)$ remains on the 
boundary of $\JJ^-(p)$ for all positive  values of $s$ or there exists a  a value
$ s_*(\om)$ beyond which  the points $\ga_\om(s)$ are no longer on the boundary
of $\JJ^{-}(p)$ but rather in its interior, see \cite{HE}.   Thus    $\NN^{-}(p)$ 
 is a smooth manifold at all points except the vertex $p$ and
 the  terminal points of its past null geodesic generators. 
Indeed,  at a terminal point  $q$  there exists  a null geodesic through 
$q$ which  fails to be  in  $\NN^-(p)$  past  $q$. This
implies  that the tangent space $T_q(\NN^-(p))$ contains the   past  tangent
 direction  of the null geodesic but not its opposite. This means   
  that $\NN^-(p)$ must be  singular at $q$. 
In what follows we shall denote by $\dot{\NN}^{-}(p)$ the regular part of $\NN^{-}(p)$,
that is the part  with its  terminal points removed. 
Clearly the null geodesic 
vectorfield $\L$ is well-defined and smooth on $\dot\NN^-(p)$.

 The parameter $s$ in the definition of  
 $\ga_\om$ is 
 an affine parameter on $\dot{\NN}^{-}(p)$, i.e.
\be{eq:affine-s}
\L(s)=1,\qquad s(p)=0.
\end{equation}
Let  $\ga$ denote the  degenerate metric induced by 
$\g$ on $\Null(p)$. Clearly  $\ga(\L,X)=0$ for any $X\in T\Null(p)$.
 Let $\chi$  denote the null second fundamental form of $\Null(p)$,
\be{eq:null-fund}
\chi(X,Y)=\g(\D_X \L, Y).
\end{equation}
where $X,Y$ are vector-fields tangent to $\Null(p)$
 and $\D$ denote the covariant derivative on $(\M,\g)$.
 Clearly $\chi$ is symmetric and 
 $\chi(\L,X)=0$ for any $X\in T\,\Null(p)$. This allows us to 
define $\trch$ as the trace of $\chi$ relative to 
  $\ga$.

Given a point
$q\in \Null(p)\setminus\{p\}$,  we can define a null  conjugate $\LLb$ to $\L$
 such that,
\be{eq:null-pair-gen}
\g(\L,\LLb)=-2,\qquad \g(\L,\L)=\g(\LLb,\LLb)=0.
\end{equation}
and further complement it by vectors 
 $(e_1,e_2)$ with the property that 
\be{eq:null-pair-gen2}
\g(\L,e_a)=\g(\LLb,e_a)=0,\qquad \g(e_a,e_b)=\de_{ab},\qquad a,b=1,2.
\end{equation}
The vectors $(\L,\LLb,e_1,e_2)$ can be locally extended to a neighborhood of a point 
$q\in \dot \NN^-(p)\setminus \{p\}$ to form a smooth local null frame.
Relative to such a frame the only non-vanishing components of the null second fundamental
form $\chi$ are $ \chi_{ab}=\g(\D_{e_a}\L, e_b)=\chi_{ba}$.
We can introduce the other  frame  coefficients,
\be{eq:null lapse-nullsecondf}
\chib_{ab}=\g(\D_{e_a}\LLb, e_b),\qquad
 \ze _a=\f12 \g( \D_a \L, \LLb),\qquad
 \etab _a=\f12 \g(e_a , \D_\L\LLb)
 \end{equation}
 Note that, in general, $\chib_{ab}$ is not symmetric.

{\bf Remark}. A  canonical way to define a  null geodesic 
 conjugate   is to take $\LLb$ the unique null vectorfield 
orthogonal to the  level  surfaces $S_s$   defined 
 by the affine parameter $s$. We refer to the
corresponding null pair as a  null geodesic pair.  We can also choose 
$e_1, e_2$  to be tangent to $S_s$.  Note that in that case
 $\chib $ is symmetric. We also note that in a neighborhood
of $p$ where $\NN^-(p)$ coincides with its regular part $\dot\NN^-(p)$
the geodesic null frame defined above is smooth away from the point $p$.

For the purpose of constructing our Kirchoff-Sobolev  parametrix we shall
make, in addition to {\bf A1}  the following assumption.

{\bf A2}.\quad {\it  We assume that  $\NN^{-}(q)$ coincides with $\dot\NN^-(q)$
past the space-like hypersurface $\Si$ for any point $q$ in a neighborhood of $p$.}

\subsection{Optical function.}\label{sect:optical} To make sense of
our Kirchoff-Sobolev formula we need  to define an optical
function\footnote{ i.e. a function which verifies \eqref{eq:Eikonal}} $u$, in a neighborhood
of $\Null(p)$,
  such that it  vanishes identically  on $\dot\NN^{-}(p)$.  
 We  define $u$ 
uniquely  relative to  the  time-like vector $\T_p$ as 
follows:

 Let  $\ep>0$ a small number and  $\Ga_\ep:(1-\ep, 1+\ep) \to \M$ denote  the
  timelike geodesic  from $p$  such that  $\Ga_\ep(1)=p$ 
 and $\Ga_\ep'(1)=\T_p$. 
From every point $q$  of $\Ga_\ep$ 
let $\NN^{-}(q)$  be the boundary of the past  set of $q$.
In view of assumption {\bf A2}  for all sufficiently small $\ep>0$,
$\NN^{-}(q)$ coincides with its regular part $\dot\NN^-(q)$ 
to the future  $\Si^+$ of  $\Si$.

We now  set $u$ to be the function, constant on each $\Null(q)$,
such that for $q=\Ga_\ep(t)$,
$$u|_{\Null(q)}=t-1.$$
This defines a smooth function 
$u$ which vanishes on $\Null(p)$  and verifies the eikonal equation \eqref{eq:Eikonal}
$$
\g^{\a\b}\pr_\a u\,\pr_\b u =0,
$$
in a neighborhood $\DD_\ep$ of $\Null(p)\cap\Si^+$. Observe that the null geodesic 
 vectorfield 
$\L=\g^{\a\b}\pr_\b u\, \pr_\a$  extends   the vectorfield  in
\eqref{eq:geod-null} to
$\DD$.  It   verifies the normalization condition,
$$\g(\L, \T_p)=\T_p(u)=1,$$
at all points of $\DD_\ep$.
 We can thus
 extend the definition \eqref{eq:null-fund}
of the null second fundamental form 
 $\chi$ and its trace $\trch$  at every point in $\DD_\ep$.

We can introduce local coordinates around any point in
$r\in\DD_\ep$  by considering the unique null geodesic
$\ga_{\om,q}$, with  $\om \in \SSS^2$,  which initiates at $q\in \Ga_\ep$ 
and passes through $r$ at value $s$ of its  affine parameter. Denoting by  $u$
 the value corresponding to the null cone $\Null(q)$ we see that
  $r$ is determined by the coordinates $u, s$ and 
 $\om\in \SSS^2$.

\subsection{Dirac measure on $\Null(p)$.} Given  our smooth optical  function $u$, 
defined in the 
neighborhood
$\DD_\ep$  of 
$\Null(p)\cap\Si^{+}$,  and a  distribution $\mu$   on the real line  $\RRR$,
supported at the origin, we  can define  the pull-back distribution
  $u^*(\mu)=\mu\circ u$ on $\DD_\ep\subset\M$ in the usual sense
of distribution theory. In the particular case when $\mu$
is  either the Dirac
 measure  $\de_0$ or its derivatives $\de_0', \de_0'',\ldots,$ we 
denote the corresponding distributions on $\M$ by 
$\de(u)$, $\de'(u)$, $\de''(u), \ldots.$  We can thus make sense of calculations
 such as,
$$\D_\a \de(u)=\de'(u)\D_\a u,\qquad \D_\a\D_\b\big(\de(u)\big)=\de''(u)\D_\a u \D_\b u + 
\de'(u) \D_\a\D_\b u$$
Clearly $\de(u), \de'(u),\ldots$ are supported  on  $\Null(p)\cap \DD_\ep$.   We can 
use the definition
of $\de(u)$ to define the integral along $\Null(p)$ of any  continuous function 
$f$ supported in $\DD_\ep$ as follows.

{\bf Definition.}
\quad Given a  continuous function $f$ 
supported in $\DD_\ep$ we  define its integral
on $\Null(p)$ by,
\be{eq:define-integral}
\int_{\Null(p)} f=<\de(u),f>
\end{equation}
\begin{proposition}\label{prop:integral-null-f}
The  definition \eqref{eq:define-integral}  
depends only   on the  restriction of $f$ 
to $\Null(p)$ and the normalization  condition \eqref{eq:norm-geod}
used in the definition of the null geodesic generator $\L$.
\end{proposition}

\begin{proof}:\quad We may assume without loss of generality
that $f$ is supported in   the domain  $\DD_\ep$, which can be
parametrized by the coordinates $u,s$ and $\om\in \SSS^2$
as described above. We can the easily calculate,
according to the definition of $\de(u)$ and coarea formula,
$$<\de(u), f>=\int_0^\infty\int_{\SSS^2} f(0,s,\om) ds d a_s$$
 where $da_s$ denotes the area element
on the $2$- surfaces $S_s$ of constant $s$.
\end{proof}

\subsection{ Kirchoff-Sobolev parametrix} Consider $\J_p$ to be a fixed k-tensor
at $p$ and let $\A$ be the unique
k -tensor-field defined along $\Null(p)$ which verifies the  linear transport equation,
\begin{equation}\label{eq:Transp}
\D_\L \A + \frac 12 \trch \A=0,\qquad (s\A)(p)=\J_p
\end{equation}
with  $s$ the affine parameter  \eqref{eq:affine-s}.
The tensor-field $\A$ can be extended  smoothly\footnote{We can in fact 
extended it  canonically by  solving the same transport  equation along
 $\Null(q)$, with  $q\in
\Ga_\ep$ and  initial data  $s\A(q)=\J_q$ where $J_q$ is an arbitrary smooth tensor-field
coinciding with
$\J_p$ at $p=q$ and $s$ the afine parameter along
$\Null(q)$. Note that, so defined, the tensor-field $\A$ is smooth away from the axis
$\Ga_\ep$.} to a small  neighborhood of $\Null(p)$. 
We can now define the distribution, or current, in
$\Si^{+}$,
\be{eq:Ade(u)}
<\A\de(u), \F>=<\de(u), \g(\A, \F)>
\end{equation}
for an arbitrary, smooth, 
$k$-tensor-field $\F$  supported in $\Si^+$. Here $\g(\A, \F)$ denotes
the full contraction of the  $k$ -tensor-fields  $\A$ and $\F$ with respect
to the space-time metric $\g$.
Observe that  the current $\A\de(u)$
 depends only on the choice of  $\T_p$ and $\J_p$
and not on the particular extensions of $u$ and $\A$.

In what follows we identify the space of $k$-tensors at $p$ and its dual with the help
of the  metric $\g$. 

{\bf Definition}. \quad {\it  We call\,\,\, $\KK^{-}_{p}=\KK^{-}_{p, \J_p}$, a $k$-tensor-field
distribution with  values in the space of $k$-tensors at $p$, defined by the formula 
$\KK^-_{p,\J_p}=\A\de(u)$,  with $\A$ defined by \eqref{eq:Transp}, 
the retarded Kirchoff-Sobolev 
parametrix at the point $p$, corresponding to $\J_p$.  If $\,\Psi$ is a solution 
 of the equation  $\square_\g\Psi=\F$, with $\F$ supported in $\Si^{+}$,
 we denote by $\Psi_{\F, \J_p}(p)$ the  $k$-tensor at $p$ defined by the integral,
\be{eq:Psi-F}
\Psi_{\F,\J_p}(p)=<\KK^{-}_{p,\J_p}, \F> =\int_{\Null(p)} \g(\A, \F).
\end{equation}} 

In the  case of the  scalar wave equation $\square_\g\psi=f$
 we can choose $\A$ to be the  scalar solution of \eqref{eq:Transp}
with initial data $(s\A)(p)=1$. In that case we have $\KK_{p}^-=\A\de(u)$
and
$$\psi_{f}(p)=<\KK^{-}_{p}, f>=\int_{\Null(p)} \A f.$$
 In the particular
case of
 Minkowski space we can easily identify $\A=\A_p(q)$ with the term $|x-y|^{-1}$
where $q=(s,y)\in \NN^{-}(p)$ and $p=(t,x)$.

\subsection{Time foliation near vertex}\label{sec:t-foliation} Returning to the construction
of $u$ in subsection \eqref{sect:optical} we observe that
the parameter $t$ along the 
geodesic $\Ga_\ep$ can be extended to a local, equidistant\footnote{ With the lapse function 
of the foliation identically one.}
 time foliation $\Si_t$, $t\in [1-\ep, 1+\ep]$
which covers  a whole neighborhood of  the point $p$,
such that $p\in \Si_1$. Indeed,  starting with  a fixed  
 spacelike hypersurface $\Si_1$ through
$p$, orthogonal to  the future unit timelike vectorfield $\T_p$,
 we can  define this geodesic foliation
using the timelike geodesics normal to $\Si_1$.  
In particular, for all $t\in [1-\ep, 1]$, if we denote by $\Om_\ep$ 
 the set 
\be{eq:D-ep}
\Om_\ep= \left (\JJ^{-}(p)\cap \Si_+\right)\setminus \cup_{t\in [1-\ep,1]}\, \Si_t
\end{equation}
then its  boundary is given by 
$$\pr \Om_\ep=\NN_\ep^{-}(p)\cup D_{1-\ep}\cup D$$
where $\NN_\ep^{-}(p)$ is the portion of $\NN^{-}(p)$
to the future of $\Si$ and the past of $\Si_{1-\ep}$, 
$D_{1-\ep}=\JJ^{-}(p)\cap \Si_{1-\ep}$ and $D=\JJ^-(p)\cap \Si$.

Let $\T=\D t$ denote the future,  unit normal to the foliation
$\Si_t$, defined in a neighborhood of $p$.  We define the  null lapse function $\varphi$ and the
second  fundamental form $k$ associated to $\Si_t$:
\be{eq:null-lapse-k}
\varphi^{-1}=\T(u)=\g(\L,\T),\quad k(X,Y)=\g(\D_X \T,Y),\quad \forall X,Y\in T\Si_t.
\end{equation}
 Clearly $\varphi(p)=1$. Since  $\T$ is a locally smooth vectorfield, $k$ is a
smooth  symmetric $2$-tensor. In particular,
$$
\|k\|_{L^\infty}\le C
$$
for some constant $C$. Similarly, since $u$ is a smooth optical function 
and $\varphi(p)=1$, the lapse $\varphi$ is a smooth bounded function 
in a neighborhood of $p$. in particular,
\be{eq:varphi-to-p}
|\varphi(q)-1|\to 0,\quad q\to p.
\end{equation}
We now recall the Raychaudhuri equation satisfied by $\trch$
along $\Null(p)$,
\be{eq:Raychadhouri}
\frac{d}{ds}(\trch) +\frac 12 (\trch)^2 =-|\chih|^2 - {\bf Ric}(\L,\L).
\end{equation}
with $s$ the afine parameter of $\L$ and $\chih$ the traceless
part of $\chi$.

The behavior of the function $\trch$ at the vertex $p$ is determined by 
the conditions  \be{eq:trch-p}
(s\trch)(p)=2,\qquad \chih(p)=0.
\end{equation} Integrating the 
Raychaudhuri equation one can easily  deduce that, 
\be{eq:trch-near-p}
|\trch(q)-\frac 2s|\to 0,\quad q\to p.
\end{equation}

Consider the time function $t$ restricted to $\NN^-(p)$. Then 
\be{eq:lapse-phi}
\frac {\partial t}{\partial s} = \L(t) = \g(\L,\T)=\varphi^{-1}
\end{equation}
The area $|S_t(p)|$ of the 2-d surfaces $S_t(p)=\Si_t\cap \NN^{-}(p)$ 
obeys the equation
$$
\frac d{dt} |S_t(p)|= \int_{S_t(p)} \varphi\, \trch\, da_\ga. 
$$
This and the behavior of $\trch$ and $\phi$ near $p$  ($t(p)=1$) imply that 
\be{eq:area-St-near-p}
|S_t(p)|=4\pi (t-1)^2 + O(|t-1|^3)
\end{equation}
On the other hand from \eqref{eq:lapse-phi} and \eqref{eq:varphi-to-p},
$t-1=s +o(s)$, which implies that 
\be{eq:area-St-near-p2}
|S_t(p)|=4\pi s^2 + o(s^2)
\end{equation}

We shall also make use of the following 
 simple variation of  proposition \ref{prop:integral-null-f}.
\begin{proposition}
Let $t$ be a  regular time function defined on $\Null(p)$ with $t(p)=1$
and equal $t_0<1$ on  $\Si_0\cap\Null(p)$, where $\Si_0$ is an arbitrary 
 spacelike hypersurface on $\Null(p)\cap\Si^+$.
 Assume that   $\varphi =\frac{dt}{ds}<0$. 
 Then, for every test function $f$,
compactly supported in $\JJ^{+}(\Si_0)$
\be{eq:integr-on-Null}
< \de(u), \psi>=\int_{t_0}^1\int_{S_t}f \varphi  dt da_t
\end{equation}
where  $S_t$ denotes the level surfaces of  $t$ and $da_t$
the corresponding area element.
\end{proposition}

\begin{proof}:\quad  The result follows easily  by first extending 
  $t$ and $u$ to  a neighborhood $\DD$ of $\Null(p)\cap\Si^+$ 
and the applying the coarea formula as above.
\end{proof}
\subsection{Statement of the result} We consider a space-like hypersurface
 $\Si\subset \M$ and a point $p\in \Si^{+}=\JJ^{+}(\Si)$ such that the assumptions
{\bf A1}-{\bf A2} are satisfied.
\begin{theorem}\label{thm:KS1}Let $\Psi$ be a solution of the equation $\square_\g\Psi=\F$
with $\F$ a $k$-tensor-field supported in $\Si^{+}$. Then for any $k$-tensor $\J_p$ at $p$,
\begin{equation}\label{eq:Hadamard} 
\Psi(p)= \Psi_{\F,\J_p}(p) + \int_{\Null(p)}\g ({\cal E},\Psi),
\end{equation}
where 
$$
 \Psi_{\F,\J_p}(p) =<\KK^-_{p,\J_p},\F>=\int_{\Null(p)} \g (\A, \F)
$$
and  $\A$ verifies \eqref{eq:Transp}. The smooth error  term ${\cal E}  $ depends only on $\J_p$, 
the  geometry of the truncated  null cone $\Null(p)\cap \Si^+\subset \M$,
and  the  ambient  spacetime 
curvature $\R$ restricted to   $\Null(p)$. 

In the particular case of a scalar wave equation $\square_\g\psi=f$, 
$\A$  and ${\cal E}$  are  scalar functions on $\Null(p)$ and,
$$
\psi(p)= \psi_f(p) + \int_{\NN^{-}(p)}{\cal E}  \psi,
$$
\end{theorem}
The precise expression of the error term $\EE$ will be given in Theorem
\ref{thm:KS2}

 \section{Derivation of Kirchoff-Sobolev formula}

\subsection{Covariant derivatives of space-time tensors}
\label{seq:cov-deriv}
As is well known there is no canonical way to define a restriction of the
space-time  covariant 
derivative $\D$ to a null hypersurface. This is due to the absence of a canonical 
projection of a tangent space $T_q{\M}$,  $q\in \Null(p)$, onto 
 the tangent space $T_q(\Null(p))$. This projection can be fixed, however, by a choice 
 of a null conjugate $\LLb$, i.e.  a null vector such that
$\g(\L,\LLb)=-2$.  With this choice we define an induced covariant derivative
$\Db$ on $\Null(p)$:
$$
\Db_X Y = \D_X Y +\frac 12 \chi(X,Y) \LLb,\qquad \forall X,Y\in T\Null(p)
$$
For example, if we choose  $X,Y$ to be the elements $(e_a)_{a=1,2}$
of a null frame $(\L,\LLb, e_1,e_2)$,
$$
\D_ae_b=\Db_a e_b+\f12 \chi_{ab}\, \LLb.
$$

We now make sense of covariant derivatives of space-time 
tensors along $\Null(p)$. 
We start by defining a covariant derivative $\Dh$ of a space-time 1-form $A_\mu$
defined on $\Null(p)$. Thus we view $A$ as a section of the vector bundle\footnote{
with the
covariant derivative  denoted  by $\D$} $T^*\M$ 
over $\Null(p)$, endowed with the induced covariant derivative $\Db$.
We interpret the covariant derivative $\Dh A$ of $A$ along $\NN^{-}(p)$ as a 1-form 
on $\Null(p)$ with values in $T^*\M$.  Thus, for every vectorfield $X\in T \Null(p)$
and any vectorfield $Z$ in $ T\M$,
\beaa
\Dh A(X; Z)=\Dh_X A(Z):=X\big(A(Z)\big)-A(\D_X Z)
\eeaa
We also write,
$$
(\Dh_X A)_\mu = X^a \D_a A_\mu,\qquad \forall X\in T\Null(p). 
$$
 We  define $\Dh^2 A$,  an $\Null(p)$ 2-tensor of second   covariant derivatives of $A$ along 
 $\Null(p)$ with values in $T^* {\M}$,
by the formula,
\beaa
\Dh^2 A(X, Y; Z)=(\Dh_X \Dh A)(Y;Z)=X(\Dh A(Y;Z))-\Dh A(\Db_X Y; Z)-
\Dh A(Y; \D_X Z)
\eeaa
or simply,
$$
\Dh^2 A_\mu (X,Y) = (\Dh_X (\Dh_Y A))_\mu - (\Dh_{\Db_X Y} A)_\mu
$$
These definitions can be easily extended to higher
covariant derivatives along $\Null(p)$ and 
to higher order tensors $A$.

\subsection{Kirchoff-Sobolev current}
 Consider the current $\KK^{-}_{p,\J_p}=\A\de(u)$ 
defined in \eqref{eq:Ade(u)}. Recall that $\J_p$
is an arbitrary $k$-tensor at $p$ and $\A$ is a $k$-tensor-field verifying the transport equation,
along $\Null(p)$,
\be{eq:transportA}
\D_\L \A+\f12 \A\trch =0,\qquad s\A(s)|_{s=0}=\J_p.
\end{equation}
Also $ \L=\g^{\mu\nu}\pr_\nu u\,\pr_\mu$, $\g^{\mu\nu}\pr_\nu u\,\pr_\mu u=0$ and $\L(s)=1$.
Let $\LLb$ be an arbitrary local null conjugate to $\L$, i.e. $\LLb(u)=\g(\L,\LLb)=-2$. 
Calculating relative to an arbitrary null frame
we easily check  that $\square_\g u= \trch.$
Formally we thus  have,
\beaa
\square_\g \big(\A\de(u)\big)&=& \square_\g \A\de(u)+
(\g^{\mu\nu}\D_\mu \A \D_\nu u+\A\square_\g u)\de'(u)+
\A(\g^{\mu\nu}\pr_\nu u\,\pr_\mu u)\de''(u)\\
&=& \square_\g \A\de(u)+
(-\LLb(u) \D_\L \A+\A\square_\g u )\de'(u)\\
&=& \square_\g \A\de(u)+2(\D_\L \A+
\f12 \A\trch)\de'(u)
\eeaa
Hence,
\be{eq:fund-sol}
 \square_\g \big(\A\de(u)\big)=\square_\g \A\de(u)+ 2(\D_\L \A+
\f12 \A\trch)\de'(u).
\end{equation}

 Observe that the above calculation
 does not depend on the choice of $\LLb$. 

Since
\be{eq:transport-sA}
\D_\L(s\A)= -\f12  (\trch-\frac{2}{s})s\A
\end{equation}
we have, in view of \eqref{eq:trch-near-p},
that along $\Null(p)$, 
\be{eq:asympt-A}
|s\A(q)-\J_p|\to 0,\qquad s\to 0.
\end{equation}

We shall next apply $\KK^{-}_{p,\J_p}=\A\,\de(u)$
to the equation $\square\psi=\F$
in the sense of distributions, 
\bea
\label{eq:Hadamard1}
\int_{\Si^+}\g (\A\,\de(u) , \square_\g \Psi)=\int_{\Si^+}\g (\square\big(\A\,\de(u)\big) , \F)
\eea
where $\Si^{+}=\JJ^{+}(\Si)$ is the future of the initial
hypersurface $\Si$. We assume that $\Psi$  has zero data on $\Si$ and that 
$\F$ is supported in $\Si^+$.
Our next goal is to integrate by parts on  the left hand side of \eqref{eq:Hadamard1}.
We first decompose $\square_\g \A $,
for an arbitrary tensor-field $\A$, relative 
to our null frame \eqref{eq:null-pair-gen} - \eqref{eq:null-pair-gen2}.
For simplicity we assume that $\A_\mu$ is a 
one tensor, the general case can be treated
 in the same  manner.  We recall the definition of the Ricci coefficients
\eqref{eq:null lapse-nullsecondf},
$$\chib_{ab}=\g(\D_{e_a}\LLb, e_b),\qquad
 \ze _a=\f12 \g( \D_a \L, \LLb),\qquad
 \etab _a=\f12 \g(e_a , \D_\L\LLb)$$
and also introduce,
\be{eq:omega}
\om=-\frac{1}{4}\g(\D_\LLb\LLb, \L)
\end{equation}
which is well defined in a neighborhood $\DD_\ep$
 of $\Null(p)$, see subsection \ref{sect:optical}.
 Using  also the notation in subsection \eqref{seq:cov-deriv}
we derive:
\beaa
\square_\g A_\mu &=&\g^{\a\b } \D^2_{\a\b} A_\mu=
 -\f12 \D^2_{\L\LLb} A_\mu -\f12 \D^2_{\LLb \L} A_\mu
+ \de^{ab} \D^2_{ab} A_\mu
\eeaa
Now, $\D_b A_\mu=\Dh_b A_\mu$ and, since 
$ \D_a e_b=\Db_a e_b +\frac 12
\chi_{ab} \LLb$,
\beaa
\D^2_{ab}A_\mu&=& e_a (\Dh_b A_\mu)-\D_{\D_a e_b} A_\mu\\
&=& e_a (\Dh_b A_\mu)-\D_{\Db_a e_b} A_\mu
- \frac 12\chi_{ab}\D_\LLb A_\mu\\
&=&\Dh^2_{ab} A_\mu - \frac 12\chi_{ab}\D_\LLb A_\mu
\eeaa
Hence, denoting $\laph A_\mu=\de^{ab}\Dh^2_{ab} A_\mu$, 
\beaa
 \de^{ab} \D^2_{ab} A_\mu&=&\laph A_\mu-\f12 \trch\D_\LLb A_\mu
\eeaa
On the other hand,
\beaa
\D^2_{\L\LLb} A_\mu&=&\D^2_{\LLb\L} A_\mu
+\R_\mu^{\,\,\la}\,_{\,\, \L \LLb}A_\la\\
&=&\D_\LLb\D_\L A_\mu- 2\ze_a\Dh_a A_\mu+2\om \Dh_\L A_\mu 
+\R_\mu^{\,\,\la}\,_{\,\, \L \LLb}A_\la\
\eeaa
Henceforth,
\bea
\square_\g A_\mu &=&-\D_\LLb\D_\L A_\mu+\laph A_\mu 
 +\ze_a\c\Dh_a A_\mu\label{eq:square-decomp}\\
&-&\om \Dh_\L
A_\mu -\f12 \trch\D_\LLb A_\mu
-\f12 \R_\mu^{\,\,\la}\,_{\,\, \L\LLb}A_\la\nn
\eea
\begin{remark}In the case when
$\A$ is a scalar  formula \eqref{eq:square-decomp} becomes,
simply,
\bea
\square_\g A_\mu &=&-\D_\LLb\D_\L A+\laph A 
 +\ze_a\c\Dh_a A
-\om \Dh_\L
A -\f12 \trch\D_\LLb A\label{eq:square-decomp-scalar}
\eea 
\label{remark:square-decomp-scalar}
\end{remark}
\subsection{Integration by parts}
In view of \eqref{eq:fund-sol}
 we have,
\bea
\square_\g \big(\A\de(u)\big)&=& \square_\g \A\,\de(u) + 
(2\Dh_\L A+\trch \A) \de'(u)\label{square-of-Ade}
\eea
where for $u=0$,
$$
\Dh_\L \A+\f12 \A\trch =0.
$$
According to 
section \ref{sec:t-foliation} 
we have defined a time foliation $\Si_t$, 
 $t\in[1-\ep,1+\ep]$,  in a neighborhood 
of the vertex $p$ with $p\in \Si_1$  such that the boundary of the set $\Om_\ep= 
\JJ^{-}(p)\setminus \cup_{t\in [1-\ep,1]}\, \Si_t$
is given by,
$$\pr \Om_\ep=\Null_\ep(p)\cup D_{1-\ep}\cup D. $$
Here $\Null_\ep(p)$ is the portion of $\Null(p)$
to the future of $\Si_{1-\ep}$ , $D_{1-\ep}=\JJ^{-}(p)\cap \Si_{1-\ep}$
and $D=\JJ^{-}(p)\cap \Si$.
As before we  denote by $\T$ the future  unit normal to the surfaces 
 $\Si_0=\Si$ and $\Si_{1-\ep}$. Note that $\Si_{1-\ep}$ needs
 only  be defined locally, for small $\ep>0$. Note also that, for $\ep=0$, $\T$
 coincides with $\T_p$ as defined in section \ref{se:null-cone}.
Clearly,
\beaa
\int_{\Om}\g (\A\,\de(u), \Box_\g \Psi)=\lim_{\ep\to 0}\int_{\Om_\ep}\g (\A\,\de(u), \Box_\g \Psi)
\eeaa
where $\Om=\JJ^{-}(p)\cap \Si_+$.
 Due to the presence of $\de(u)$ and the fact that $\psi$ is
supported in $\Si_+$ 
we may  assume in what follows  that  all functions 
we deal with in the calculation below  are supported 
in the set $\Om=\JJ^{-}(p)\cap J^+(\Si)$. Thus  the boundary of the
intersection of their supports with $\Om_\ep$ is included in $\Si_{1-\eps}$.

\begin{lemma}
Let $ F, G$ be two tensor-fields of the same rank and $F$ is a distribution supported in $\Omega$. 
Then
\beaa\int_{\Om_\ep} \g (F,\Box_\g  G)=\int_{\Om_\ep} \g(\Box_\g  F,G)-
\int_{D_t} \big(  \g(F,\D_\T G)-\g (G,\D_\T F)\big)\bigg|_{t=0}^{t=1-\ep},
 \eeaa
 where  $D_0=D$.
\end{lemma}
\begin{proof}:\quad Indeed,
\beaa
\g( F,\Box_\g  G) -\g (\Box_\g F,G)=\D^\a \g(F, \D_\a G)- \D^\a\g (\D_\a F, G)
\eeaa
Thus,
\beaa
\int_{\Om_\ep} \left(\g(F, \Box_\g G)  -\g(\Box_\g  F,G)\right )&=&
\int_{\Om_\ep}\D^\a\big( \g(F, \D_\a G)- \g(\D_\a F,G)\big)\\
&=&-\int_{D_t} \T^\a \big( \g(F, \D_\a G)-\g(\D_\a F, G)\big) \bigg|^{1-\ep}_0
\eeaa
\end{proof}
We  now write, 
\beaa
\int_{\Om_\ep}\g\left(\A\,\de(u),\Box_\g \Psi\right)&=&\int_{\Om_\ep}
\g\left(\Box_\g \big(\A\,\de(u)\big), \Psi\right)\\
&-& \int_{D_t}\g(\A\de(u), \D_\T \Psi)\bigg|_0^{1-\ep}+
\int_{D_t}\g\left (\D_\T\big(\A\de(u)\big),
\Psi\right)\bigg|_0^{1-\ep}\\ &=&\int_{\Om_\ep}\g\left(\Box_\g\big (\A\,\de(u)\big),  
\Psi\right)+ I_\ep+J_\ep, 
\eeaa
where $I_\ep$ and $J_\ep$  denote the boundary terms 
 on $D_{1-\ep}$. The  term corresponding to  $D_0$ vanishes due
 to the zero data assumption for $\Psi$.  
\begin{proposition}
We have 
\beaa
   I_\ep\to 0,\qquad   J_\ep \to -
4\pi \g(\Psi(p), \J_p) \quad \mbox{as }\quad \ep\to 0.
\eeaa
Thus,
\bea
\int_{\Om}\g\left(\A\,\de(u),\Box_\g \Psi\right)&=&
\lim_{\ep\to 0}\int_{\Om_\ep}\g\left(\Box_\g\big (\A\,\de(u)\big), \Psi\right) -
4\pi\g\left( \Psi(p), \J_p\right)
\label{eq:integr-parts-form} 
\eea
\end{proposition}
\begin{proof}:\quad 
We   analyze the  boundary terms  $I_\ep, J_\ep.$
Clearly,
\beaa
I_\ep= -\int_{D_{1-\ep}}\g\left(\A\de(u), \D_\T \Psi\right)&=&
-\int_{\Null(p)\cap D_{1-\eps}}\g(\A, \D_\T \Psi )\,\varphi  \,da_\ga
\eeaa
where,  see \eqref{eq:null-lapse-k},
 $\varphi=|\D_\T u|^{-1}$ is the null lapse 
 and $da_{\ga}$ is the area element 
of the 2-surface $S_{1-\ep}(p)=\Null(p)\cap D_{1-\eps}$.  Recall that 
according to  \eqref{eq:area-St-near-p} the area
$|S_{1-\ep}(p)|\les \ep^2$.

Now,
\beaa
|I_\ep|& \les &\|\varphi\|_{L^\infty}\big(\int_{S_{1-\ep}}|\A|^2 da_\ga\big)^{1/2} 
\big(\int_{S_{1-\ep}}|\D_\T\Psi|^2 da_\ga\big)^{1/2} \\
&\les& 
\|\varphi\|_{L^\infty}\|\D_\T\Psi\|_{L^\infty}\|\A\|_{L^2(S_{1-\ep}(p))}|S_{1-\ep}(p)|^{1/2}\\
&\les& \ep\, \|\varphi\|_{L^\infty}\|\D_\T\Psi\|_{L^\infty}\|\A\|_{L^2(S_{1-\ep}(p))}
\eeaa
Recalling \eqref{eq:asympt-A} and \eqref{eq:area-St-near-p} we easily see that 
$\|\A\|_{L^2(S_{1-\ep})}$ is bounded as $\ep\to 0$.
Thus, for a smooth  tensor-field $\Psi$, we clearly
have,
$$I_\ep \to 0, \qquad \mbox{as}\quad \ep \to 0.$$

We now consider the second boundary term,
\beaa
J_\ep&=&\int_{D_{1-\ep}}\g\left(\D_\T\big(\A\de(u)\big),\Psi\right)\\
&=&\int_{D_{1-\ep}}\de(u)\, \g\left( \D_\T\A,\Psi\right) +\int_{D_{1-\ep}}\de'(u) \D_\T u\, \,\g(\A, \Psi)\\
&=&\int_{D_{1-\ep}}\de(u) \,\g\left(\D_\T\A, \Psi\right) +\int_{D_{1-\ep}}\de'(u)\, \varphi^{-1} \,\g(\A, \Psi)\\
&=&J_\ep^1+J_\ep^2
\eeaa
 If $N=\varphi \L+\T$ denotes the unit 
 normal\footnote{A priori,  the vectorfield $N$ is defined only  on 
$\Null(p)\cap\Si_{1-\eps}$, it can however be extended locally as a unit 
normal to the foliation of 2-d surfaces 
$\{u=const\}\cap \Si_{1-\eps}$.  } to $S_{1-\ep}=\Null(p)\cap\Si_{1-\eps}$ 
 in $\Si_{1-\eps}$, then 
 $\D_N \de(u)=\de'(u) \D_N u$ and $\D_N u =\D_\T u=\varphi^{-1}$.
Hence,
\beaa
J_\ep^2=\int_{\Si_{1-\ep}}\de'(u) \varphi^{-1} \,\g(\A, \Psi)&=&
 \int_{\Si_{1-\ep}}\D_N \de(u)  \,\g(\A, \Psi)
\eeaa
We next record  the following integration by parts formulae.
\begin{lemma}
Let $X$ be a vectorfield tangent to the hyperplane $\Si_t$ and let $f, g$ be 
two scalar functions on $\Si_t$.
Denote by $\nab$ the covariant derivative  restricted 
 to $\Si_t$.  
Then,  
\begin{equation}
\int_{\Si_t} f X(g) = - \int_{\Si_t} (X(f) + \div X f) g.
\label{bypartX}
\end{equation}
In particular,
\begin{equation}
\int_{\Si_t} f N(g) = - \int_{\Si_t} (N(f) + \tr\th\, f) g,
\label{bypartN}
\end{equation} 
where $\tr\,\th$ is the mean curvature of the 2-d surfaces $S_{t,u}=
\{u=const\}\cap\Si_t$.
\end{lemma}
\begin{proof}:\quad Formula \eqref{bypartX} is  standard.  To prove 
\eqref{bypartN} observe  $\div N=\g(\nab_ N N,
N)+\sum_a\g(\nab_a N, e_a) =\tr\th$ where $\th$ is the second
 fundamental form of  the surfaces
$S_{t,u}\subset
\Si_t$.
\end{proof}
Using the lemma
we infer that,
\beaa
J_\ep^2&=& -\int_{\Si_{1-\ep}} \de(u)\big(  N \g(\,\A,\Psi) +\tr\th \,\g(\A,\Psi)\big)\\
&=& - \int_{\Si_{1-\ep}} \de(u)\big(  \g(\D_N\A,\Psi)  +\tr\th\, \g(\A,\Psi)\big)
-\int_{\Si_{1-\ep}} \de(u)\,\g(\A, \D_N \Psi)
\eeaa
Now, proceeding as for $I_\ep$, it is easy to check that 
$\int_{\Si_{1-\ep}} \de(u)\,\g(\A, \D_N \Psi)\to 0$ as $\ep\to 0$.
Hence,
\beaa
\lim_{\ep\to 0}J_\ep^2&=&- \lim_{\ep\to 0}
\int_{\Si_{1-\ep}} \de(u) \big(  \g(\D_N\A,\Psi)  +\tr\th\, \g(\A,\Psi)\big)
\eeaa
Or,
\beaa
\lim_{\ep\to 0}J_\ep&=&- \lim_{\ep\to 0}\int_{\Si_{1-\ep}} \de(u)\big(\g( \D_{N-\T} \A,\Psi)+
\tr\th\,\g(\A,\Psi)\big) 
\eeaa
Now observe  that $\L=\varphi^{-1}(N-\T)$. Hence,
$\D_{\L}\A =\varphi^{-1} \D_{N-\T}\A$.
Since  $\D_\L \A+\f12 \A\trch =0$ we infer,
$$\D_{N-\T}\A =\varphi\, \D_\L \A =-\f12 \varphi  \A\trch.$$
Hence,
\beaa
\lim_{\ep\to 0}J_\ep&=-&\lim_{\ep\to 0} \int_{\Si_{1-\ep}} \de(u)\big(-\f12 \varphi \trch+\tr\th\big)
\g(\A, \Psi)
\eeaa
On the other hand $\th_{ab}=g(\nab_a N, e_b)=\g(\D_a(\varphi \L+\T), e_b)=\varphi \chi_{ab}+k_{ab}$.
Therefore, $\tr\th=\varphi\trch + \de^{ab} k_{ab}$ and we  deduce,
\beaa
\lim_{\ep\to 0}J_\ep&=&
-\f12 \lim_{\ep\to 0} \int_{\Si_{1-\ep}}\varphi \trch  \de(u) \,\g(\A, \Psi) -
\lim_{\ep\to 0} \int_{\Si_{1-\ep}}(\de^{ab} k_{ab})\de(u)\,\g(\A, \Psi)
\eeaa
It is easy to see that the second term
of the right hand side converges to zero
for $\ep\to 0$. Indeed,
\beaa
|\int_{\Si_{1-\ep}}(\de^{ab} k_{ab})\de(u)\,\g(\A,\Psi)\, |&=&\big|\int_{S_{1-\ep}(p)} 
 (\de^{ab} k_{ab})\,\g(\A, \Psi)\, \varphi \,da_\ga\big|\\
&\les&|S_{1-\ep}(p)|^{\frac 12}\|k\|_{L^\infty}
\|\Psi\|_{L^\infty}\|\varphi|_{L^\infty}\|\A\|_{L^2(S_{1-\ep}(p))} 
\eeaa
Therefore,
\beaa
\lim_{\ep\to 0}J_\ep&=&-\f12 \lim_{\ep\to 0} \int_{S_{1-\ep}(p)}\varphi^2\, \trch\, \g( \A, \Psi)\, da_\ga
\eeaa
It is easy to check that,
\beaa
\f12 \lim_{\ep\to 0} \int_{S_{1-\ep}(p)}\varphi^2\, \trch\,\g\left(\A, (\Psi-\Psi(p)\right) da_\ga=0
\eeaa 
Therefore,
\beaa
\lim_{\ep\to 0}J_\ep&=& -\f12\lim_{\ep\to 0} \int_{S_{1-\ep}(p)}\varphi^2 \trch\, \g\left(\A,\Psi(p)\right)\,
d\a_\ga 
\eeaa
Or, since $\sup_{S_{1-\ep}(p)}|\varphi-1|\to 0$,
 and $\sup_{S_{1-\ep}(p)}|\trch-\frac{2}{s}|\to 0$
 as $\ep\to 0$, we infer 
that
\beaa
\lim_{\ep\to 0}J_\ep&=&-\ep^{-2} \lim_{\ep\to 0}\int_{S_{1-\ep}(p)} \g\left(\Psi(p), ( s\A)\right)=
- 4\pi  \lim_{r\to 0} \g\left(\Psi(p), (r\A)(r)\right)
\eeaa
Thus, using the initial condition $\lim_{s\to 0} s \A(s) = \J_p$, 
we obtain
\beaa
\lim_{\ep\to 0}J_\ep&=&-
4\pi \g(\Psi(p), \J_p)
\eeaa

\end{proof}
We now analyze the term $\int_{\Om_\eps}\g\left(\Box_\g \big(\A\,\de(u)\big), \Psi\right)$
 on the right hand side of \eqref{eq:integr-parts-form}.
In view of \eqref{square-of-Ade} we have, 
$$
\int_{\Om_\eps} \g\left(\Box_\g \big(\A\,\de(u)\big),\Psi\right)=
\int_{\Om_\eps}\de(u)\,\g\left((\Box_\g \A), \Psi\right) + 
\int_{\Om_\eps} \de'(u)\, \g\left((2\D_\L A+\trch \A), \Psi\right).
$$
Given the normalization $\LLb (u)=-2$  we have
$\de'(u)= -\frac 12\D_\LLb \de(u)$. Integrating by parts we obtain
\begin{align*}
\int_{\Om_\eps} \de'(u)\,\g\left ((2\D_\L A+\trch \A), \Psi\right) & = 
\frac 12\int_{\Om_\eps}\de(u)\, \g\left(\D_\LLb (2\D_\L A+\trch \A), \Psi\right)\\ &+
\frac 12 \int_{\Om_\eps}\de(u)\,\g\left ((2\D_\L A+\trch \A), \big (\D_\LLb \Psi
+\D^\a \LLb_\a \Psi\big )\right) \\ &+ \int_{D_t}\de(u)\,\g(\left((2\D_\L A+\trch \A),\Psi\right)
\,\g(\LLb,\T)\, |_{t=0}^{t=1-\eps}
\end{align*}
Recall that $2\D_\L A+\trch \A=0$ on the surface $u=0$. 
Therefore
the last two terms vanish and we derive,
$$
\int_{\Om_\eps}\de'(u)\,\g\left((2\D_\L A+\trch \A), \Psi\right)=
\int_{\Om_\eps} \de(u)\,\g\left(\D_\LLb(\D_\L A+\frac 12\trch \A),\Psi\right)
$$
Therefore,
\beaa
\int_{\Om_\eps} \g\left(\Box_\g \big(\A\,\de(u)\big), \Psi\right)=
\int_{\Om_\eps}\de(u)\,\g\left(\big(\Box_\g \A  + \D_\LLb(\D_\L A+\frac 12\trch \A) 
  \big), \Psi\right)
\eeaa
We now  recall \eqref{eq:square-decomp},
\beaa
\Box_\g \A &=&-\D_\LLb\D_\L \A+\laph \A 
 +\ze_a\c\Dh_a \A\\
&-&\om \Dh_\L
\A -\f12 \trch\D_\LLb \A
-\f12 \R(\c,\A, \L, \LLb)\nn
\eeaa
Therefore,
\begin{align*}
\Box_\g \A  + \D_\LLb(\D_\L A+\frac 12\trch \A)&=\laph \A +
\ze_a \Dh_a \A - \frac 12 \trch \D_\LLb \A\\
&- \om \Dh_\L \A 
+\f12 (\D_\LLb\trch) \A-\f12 \R(\c,\A, \L, \LLb)\,
\end{align*}
while since $\D_\L \A +\frac 12 \trch \A=0$ on $\Null(p)$,
$$
 -\frac 12 \trchb \D_\L \A -\om
\D_\L \A +\f12 (\D_\LLb\trch) \A =\f12(\D_\LLb\trch+\f12 \trch\trchb+2\om \trch)\A
$$
Hence, we have proved the following,
\begin{proposition}In the case of a one form   $\A$ verifying \eqref{eq:transportA},
\beaa
\Box_\g \A  + \D_\LLb(\D_\L \A+\frac 12\trch \A)&=&\De \A +
\ze_a \D_a \A\\
 &+&\f12(\D_\LLb\trch+\f12 \trch\trchb+2\om \trch)\A-\f12
\R(\c,\A, \L, \LLb),
\eeaa
where 
$$
\lap \A = \laph \A + \frac 12 \trch \Dh_\L \A=e_a (\D_a A)-\g(\D_a e_a,e_b) \D_b A
$$
is a  Laplace-Beltrami  type operator which coincides with the standard surface 
Laplace-Beltrami operator in the case when  the frame $\{e_a\}_{a=1,2}$ spans a tangent 
space of a 2-dimensional surface\footnote{As in the case of the 
geodesic foliation.}. 

In the  scalar case we have instead,  see remark \ref{remark:square-decomp-scalar},
\beaa
\Box_\g \A  + \D_\LLb(\D_\L \A+\frac 12\trch \A)&=&\De \A +
\ze_a \D_a \A\\
 &+&\f12(\D_\LLb\trch+\f12 \trch\trchb+2\om \trch)\A,
\eeaa

\end{proposition}
Using the above proposition we infer that, 
\beaa
\int_{\Om_\ep}\g\left(\Box_\g \big(\A\,\de(u)\big),\Psi\right)&=&
\int_{\Om_\ep} \de(u) \g\left(\big(\lap \A +\ze_a \D_a \A \big), \Psi\right)\\
&+&\f12\int_{\Om_\ep}\de(u)\,\big( \D_\LLb\trch+\f12 \trch\trchb+2\om \trch\big)\,\g(\A,\Psi)\\
&+&\f12 \int_{\Om_\ep} \de(u) \R(\Psi,\A, \LLb\, , \L)
\eeaa
We now make use of the following,
\begin{proposition}
Introduce the mass aspect function as in (13.1.10b) of \cite{C-K},
\be{eq:mass-asp}\mu= \D_\LLb\trch+\f12 \trch\trchb+2\om \trch
\end{equation} 
The following formula holds true relative to the standard geodesic foliation on $\Null(p)$, 
\be{eq:def-mu}
\mu=2\div\ze -\chih\c\chibh+2|\ze|^2
+\R_{\LLb\L} + \frac 12 \R(\LLb,\L,\LLb,\L)
\end{equation}
\end{proposition}
\begin{proof}:\quad See \cite{C-K}. 
\end{proof}
\begin{remark}
Note that according to \eqref{eq:def-mu} the mass aspect function $\mu$ depends 
only on the null hypersurface $\Null(p)$ and the ambient curvature $\R$.
\end{remark}
We have therefore
proved the following precise version of theorem \ref{thm:KS1},
\begin{theorem}
\label{thm:KS2}
Let $\A$ be a vectorfield verifying,
 $$
\D_\L \A+\f12 \A \trch =0,\quad  s\A(p)=\J_p \qquad{\text on}\quad u=0
$$
where $\J_p$ is a fixed vector at $p$. Then solution $\Psi$ of an inhomogeneous
vector equation $\Box_\g\Psi=\F$, in a globally hyperbolic spacetime $(\M,\g)$ satisfying 
${\bf A1}, {\bf A2}$, with zero initial data on a Cauchy hypersurface $\Si$ can be represented
by the following formula at point $p$ with $\Omega=\JJ^-(p)\cap \JJ^+(\Si)$,
\bea
4\pi \g(\Psi(p), \J_p)&=&-\int_{\Om}\de(u)\, \g(\A,\F) \label{eq:main-ident}\\
&-&\frac 12\int_{\Om}\de(u)\, \R(\Psi,\A, \LLb\, , \L)
+\int_{\Om} \de(u)\,\g\left(\big(\lap \A
+\ze_a \D_a \A \big),\Psi\right) \nn\\ &+&\f12\int_{\Om}
\de(u)\,\mu\,\g(\A,\Psi)\nn
\eea
where,
\beaa
 \R(\Psi,\A, \LLb\, , \L)=\R_{\a\b \ga\de} \LLb^\ga  \L^\de\Psi^\a\A^\b,
\eeaa
with $\R_{\a\b\ga\de}$ the components of the curvature tensor
$\R$ relative to an arbitrary frame.
\end{theorem}
\begin{remark}
Theorem \ref{thm:KS2} implies that the error term ${\cal E}$ in \eqref{eq:Hadamard} 
has the following representation
$${\cal E}_a = -\f12  \R_{a\la \ga\de}  \A^\la \LLb^\ga \L^\de+(\lap \A
+\ze_a \D_a \A )_\a + \frac \mu 2 \A_\a
$$
in the case of a vectorial wave equation. For the scalar wave equation 
$$
{\cal E} = (\lap \A
+\ze_a \D_a \A ) + \frac \mu 2 \A
$$
\end{remark}
\begin{remark}Formula \eqref{eq:main-ident} can easily be generalized to
higher order tensor wave equations. Indeed if both $\Psi$ and $\F$ are tensor-fields
of order $k$ then $\J_p$,  $\A$  and $\EE$ are also of order $k$ and,
\bea
\g(\EE,\Psi)&=&\g\big( (\lap \A
+\ze_a \D_a \A ),\Psi\big) + \frac \mu 2 \g(\A,\Psi)\\
&+&\R(\c,\c, \LLb\, , \L)\#\Psi\#\A\nn
\eea
where the last term denotes a scalar
contraction of $\R(\c,\c, \LLb\, , \L)$  with $\Psi$ and $\A$.
\end{remark}

\section{Wave equation for sections of vector bundles 
and applications to the Yang-Mills equations} 
\def\V{{\bf V}}
\def\VV{{\cal V}}
Now let $\V$ be a vector bundle over $(\M,\g)$ with a 
positive definite scalar product $<,>$ and 
a compatible connection $\la$. We may assume that $\V$ is a vector bundle  
associated to a principal bundle $P$ so that $\V=P\times_G E$ with $G$ a compact Lie group
and a vector space $E$. Let ${\cal G}$ denote the Lie algebra of $G$. 
The connection $\la$ is a ${\cal G}$ valued 1-form on $V$, which, 
locally can be viewed  as a ${\cal G}$ valued 1-form on $\M$.

We define the gauge wave operator $\Box_\g^{(\la)}$ for
sections $\Psi: \M\to\V$
$$
\Box_\g^{(\la)} \Psi= \g^{\mu\nu} \Ds_\mu \Ds_\nu \Psi,
$$
where   $\Ds_\mu = \D_\mu + [\la_\mu,\,\c\,\,]$ denotes the gauge
covariant derivative.
We  denote by $\La$  the curvature of the connection, i.e.  
 the ${\cal G}$ valued 2-form on $\M$,
$$
\La_{\a\b}=\pa_\a \la_\b-\pa_\b \la_\a + [\la_\a,\la_\b].
$$

  As before we construct  a Kirchoff-Sobolev parametrix ${\cal K}_p^-$ for
$\Box_\g^{(\la)}$ by defining $\KK_p^-=\KK_{p,\J_p}^-=\A\, \de(u)$, where $\A$ is a
section of $\V$ which verifies the  covariant transport equation\footnote{Note that here 
transport of $\A$ along integral curves of $\L$ is modulated by the action of the gauge potential $\la$.} 
\be{eq:cov-transort-A}
\Ds_\L \A +\f12 \trch \A=0
\end{equation}
with initial data $(s\A)|_{s=0}=\J_p$ and $\J_p$ is a fixed element 
of the fiber $\V_p$.
As before  we assume that  $(\M,\g)$
is  globally hyperbolic and  satisfies  
${\bf A1}, {\bf A2}$. We also assume 
that  $u$ is a solution of the eikonal equation
$\g^{\a\b}\pa_\a u\,
\pa_\b u=0$ with 
$u$, vanishing on the boundary $\NN^{-}(p)$ of the past
of $p$ in $\M$.
Repeating our calculations of section 3, we obtain the following 
analog\footnote{The only new term in the formula
 \eqref{eq:main-id} below  is due to  the commutator 
 $(\D^{(\la)}_L \D^{(\la)}_\Lb-\D^{(\la)}_\Lb \D^{(\la)}_\L)\A=[\La_{L\Lb}, \A]$
 in deriving formula  \eqref{eq:square-decomp}} of Theorem \ref{thm:KS2}. 
\begin{theorem}
\label{thm:KS3}
Let $\A$ be a section of a vector bundle $\V$ over $(\M,\g)$ verifying,
 $$
\Ds_\L \A+\f12 \A \trch =0,\quad  s\A(p)=\J_p \qquad{\text on}\quad \NN^{-}(p)
$$
where
$\J_p$ is a fixed element of $\V_p$. The solution $\Psi$ of the  inhomogeneous
gauge equation $\Box^{(\la)}_\g\Psi=\F$,  with zero initial data on a Cauchy
hypersurface $\Si$ can be represented by the following formula at point $p$ with
$\Omega=\JJ^-(p)\cap \JJ^+(\Si)$,
\bea
4\pi <\Psi(p), \J_p>&=&-\int_{\Om}\de(u)\, <\A,\F> \label{eq:main-id}\\
&-&\frac 12\int_{\Om}\de(u)\, <[\La_{\LLb \L}, \A], \Psi>
+\int_{\Om} \de(u)\,\left<\big(\lap^{(\la)} \A
+\ze_a \Ds_a^{(\la)} \A \big),\Psi\right> \nn\\ &+&\f12\int_{\Om}
\de(u)\,\mu\,\,<\A,\Psi>\nn
\eea
In particular, in Minkowski space, for general
initial data, we have the following representation:
 \bea
4\pi <\Psi(p), \J_p>&=&-\int_{\Om}\de(u)\, <\A,F> 
-\frac 12\int_{\Om}\de(u)\, <[\La_{\LLb \L}, \A], \Psi>
\label{eq:main-idM}\\&+&\int_{\Om} \de(u)\,<\lap^{(\la)} \A,\Psi>+ 
\int_{\Si} \left (<\A \de(u), \D_\T \Psi>-<\D_\T (\A\de(u)),\Psi>\right) \nn\eea
The last term represents contribution of the initial data on $\Si$.
\end{theorem}
\begin{remark}
The formula \eqref{eq:main-id} can be naturally extended to consider  
sections of the bundle $ TM\otimes\ldots\otimes TM\otimes V$. In this case
terms 
 of the form  $\int_{\Om} \de(u)\,<\R(\c,\c,\L,\LLb) \A,\Psi>$, where
$\R$ is the Riemann curvature tensor of $\g$, need to be added. The  corresponding extension of
\eqref{eq:main-idM}  does
not therefore  introduce
any additional terms.
\end{remark}
\subsection{Yang-Mills equations}
We now assume that $\la$ is a Yang-Mills connection on a 4-dimensional Lorentzian
manifold $(\M,\g)$, i.e., it verifies the equations
\begin{equation}\label{eq:Yang}
\Ds^\a \La_{\a\b}=0.
\end{equation}
The Yang-Mills equations are hyperbolic in nature and admit a Cauchy formulation,
in which the  connection $\la$ is prescribed on a Cauchy hypersurface\footnote{This
requires space-time $\M$ to be globally hyperbolic.} $\Si$ and then extended 
 as a solution of the problem \eqref{eq:Yang}. The
uniqueness and global existence for the Yang-Mills equations with smooth initial data
in 4-dimensional Minkowski space-time was established by Eardley-Moncrief, \cite{EM1},
\cite{EM2}. This result was later extended to the Yang-Mills equations on a smooth
4-dimensional globally  hyperbolic Lorentzian space-time by Chru\'sciel-Shatah,
\cite{CS}. A different proof in Minkowski space, allowing for  initial data
 with only finite energy , was
given by Klainerman-Machedon,
\cite{KM}. 

All of the above approaches were manifestly non-covariant; as they required 
 a choice of a gauge condition  for the connection $\la$.
The approach of Eardley-Moncrief was based on the fundamental solution for 
a scalar wave equation in Minkowski space (Kirchoff formula) and made use of Cronstr\"om
gauges: For any point $p$ the  connection $\la$ can be chosen to satisfy the condition
$$
(p-q)^\a\, \la_\a(q)=0
$$
The work of Chru\'sicel-Shatah relied on the Friedlander's representation of the fundamental
solution of a scalar wave equation in a curved space-time and a local analog of the Cronstr\"om
gauge. Finally, Klainerman-Machedon's proof was based on a Fourier representation of the
fundamental solution of a scalar wave equation in Minkowski space, bilinear estimates and 
the use of the Coulomb gauge:
$$
\pa^i \la_i=0
$$ 
Below we present a new simple {\it gauge independent} proof of the global existence and
uniqueness result for the $3+1$-dimensional Yang-Mills equations. The main new
ingredient is the use of a gauge covariant first order Kirchoff-Sobolev parametrix
described in Theorem \ref{thm:KS3}.
 
Differentiating the Bianchi identities $\Ds_{[\a} \La_{\b\si]}=0$ and using the
equations  we infer that the curvature $\La$ is a solution of a covariant gauge wave
equation
$$
\Box_\g^{(\la)} \La_{\a\b} = 2[\La^\si_{\,\,\a},\La_{\si\b}]+ 2 \R_{\si\a\ga\b} \La^{\si\ga} + 
\R_{\a\si} \La_{\b}^{\,\,\,\si} + \R_{\b\si} \La_{\,\,\b}^{\si}
$$
For simplicity we consider the problem in Minkowski space, although our results 
 can easily be extended  to the  general case of a globally hyperbolic smooth
Lorentzian manifold, as in \cite{CS}. The equations then  simplify,
\begin{equation}\label{eq:Mills}
\Box^{(\la)} \La_{\a\b} = 2[\La^\si_{\,\,\a}\,,\,\La_{\si\b}]
\end{equation}
Recall that the curvature $\La$ can be decomposed into its electric and magnetic 
parts $E_i=\La_{0i}$ and $H_i=\dual \La_{0i}$.  We also recall that
the total energy
$$
{\cal E}_0=\int_{\Si_t} \left (|E|^2+|H|^2\right).
$$
is conserved. Moreover, adapting the energy identity
to to the past $\JJ^{-}(p)$ of a point $p\in \Si^{+}$ we also get a 
bound on the flux of energy along $\NN^{-}(p)$. More precisely we derive,
$$
{\cal F}_p^-\le {\cal E}_0.
$$ 
where the  backward null energy flux ${\cal F}_p^-$ is  defined with the help of a
null frame $(\L,\LLb, e_a)$ centered at $p$. Without loss of generality we may assume 
that $p=(t,0)$ and denote $r=|y|$. Then  $\L=\pa_r-\pa_t$, $\LLb=-\pa_t-\pa_r$
and $e_a$ is a frame on a standard sphere ${\Bbb S}_r$. With these notations
$$
{\cal F}_p^-=\int_{\NN^-(p)} \left (|\La_{\L\LLb}|^2 + \sum_{a=1}^2 |\La_{\L a}|^2\right)
$$
As in the original approach of Eardley-Moncrief the key element of the proof of global
existence is a pointwise bound on curvature $\La$. Once this bound is established 
the remaining steps concerning  existence, propagation of regularity and uniqueness are
very standard and will be omitted. The precise statement
 concerning an $L^\infty$ bound
on $\La$ is as follows:
\begin{lemma}\label{lem:It}
There exists $\tau_*>0$ dependent only on ${\cal E}_0$ such that for any point 
$p=(t,0)$ we have 
$$
|\La(p)|\le C_{t-\tau_*},
$$
where the constant $C_{t-\tau_*}$ depends only on the solution $\La$ on a hypersurface 
$\Si_{t-\tau_*}$.
\end{lemma}
\begin{remark}
Iterations of Lemma \ref{lem:It} leads to a pointwise bound on the curvature $\La$ 
in terms of the initial data.
\end{remark}
\begin{proof}:\quad 
We fix $\tau_*>0$, whose is to be determined later, and apply the representation 
formula \eqref{eq:main-idM} in the domain $\Om=\JJ^-(p)\cap \JJ^+(\Si_{t-\tau_*})$,
\bea
4\pi <\La(p), \J>&=&-2 \int_{\Om}\de(u)\, <\A,[\La,\La]>
-\frac 12\int_{\Om}\de(u)\, <[\La_{\LLb \L},\A], \La>\nn \\
&+&\int_{\Om} \de(u)\,<\lap^{(\la)} \A,\La> \nn\\
&+& \int_{\Si_{t-\tau_*}} \left (<\A \de(u), \D_\T \La>-<\D_\T (\A\de(u)),\La>\right)
\label{eq:rep-M}
\eea
Here $\J$ is an arbitrary ${\cal G}$ valued 
anti-symmetric $2$-tensor on ${\Bbb R}^{3+1}$, $\A$ is a ${\cal G}$ valued 2-form 
on ${\Bbb R}^{3+1}$ verifying the equation\footnote{Recall that $\L=\pa_r-\pa_t$,
$\Ds_\L=\pa_r-\pa_t + [\la_\L,\c]$ and $r=|y|$.}  
$$
\Ds_\L \A + r^{-1} \A=0, \qquad (r\A)|_{r=0}=\J
$$ 
and $<,>$ denotes a positive definite  scalar product on $\La^2({\Bbb R}^{3+1})\otimes
{\cal G}$.   The last term in \eqref{eq:rep-M} depends only on the solution $\La$ on
$\Si_{t-\tau_*}$ and therefore is consistent with the claim of Lemma \ref{lem:It}.
We now observe that for $a,b\in \La^2({\Bbb R}^{3+1})\otimes {\cal G}$ 
we have 
$$
|<a,b>|\le |a|\, |b|,
$$
where $|a|$ denotes the absolute value of an element in $\La^2({\Bbb R}^{4})\otimes
{\cal G}$ relative to the  positive definite scalar product\footnote{\,${\Bbb R}^4$ here
stands for the  Euclidean 4-dimensional space.}. In what follows all the norms will be
understood to involve  the absolute value $|\c|$ on  $\La^2({\Bbb R}^{4})\otimes {\cal
G}$. 
We denote by  $\NN^-_{\tau_*}(p)$  the null boundary of $\Omega$
 to the future of $\Si_{\tau_*}$. Then 
\begin{align*}
&|\int_{\Om}\de(u)\, <\A,[\La,\La]>|\le \|r \A\|_{L^\infty(\NN^-_{\tau_*}(p))} \|r^{-1} [\La,\La]\|_{L^1(\NN^-_{\tau_*}(p))},\\
&|\int_{\Om}\de(u)\, <[\La_{\LLb \L},\A], \La>|\le  \|r
\A\|_{L^\infty(\NN^-_{\tau_*}(p))} 
\|r^{-1}\La_{\L\LLb}\|_{L^1(\NN^-_{\tau_*}(p))}
\|\La\|_{L^\infty(\NN^-_{\tau_*}(p))},\\
&|\int_{\Om} \de(u)\,<\lap^{(\la)} \A,\La>|\le  
\|\La\|_{L^\infty(\NN^-_{\tau_*}(p))}   \|\lap^{(\la)} \A\|_{L^1(\NN^-_{\tau_*}(p))} 
\end{align*}
It is easy to see, see e.g. \cite{EM2} that 
$$
|[\La,\La]|\le |\La| \left (|\La_{\L\LLb}|+\sum_{a=1}^2|\La_{\L a}|\right)
$$
and therefore
$$
 \|r^{-1} [\La,\La]\|_{L^1(\NN^-_{\tau_*}(p))}\le \tau_*^{\frac 12} 
 \|\La\|_{L^\infty(\NN^-_{\tau_*}(p))} \left ({\cal F}_p^-\right )^{\frac 12}
$$
Similarly,
$$
\|r^{-1}\La_{\L\LLb}\|_{L^1(\NN^-_{\tau_*}(p))}\le \tau_*^{\frac 12} 
\left ({\cal F}_p^-\right )^{\frac 12}
$$
To prove Lemma \ref{lem:It} it would be sufficient to show that 
\begin{equation}\label{eq:A}
\|r \A\|_{L^\infty(\NN^-_{\tau_*}(p))}\les |\J|,\qquad 
 \|\lap^\la \A\|_{L^1(\NN^-_{\tau_*}(p))}\les |\J| \, \left (\tau_*^{\frac 32} 
\left ({\cal F}_p^-\right )^{\frac 12} +\tau_* {\cal F}_p^-\right)
\end{equation}
and then  choose $\tau_*<<(1+{\cal E}_0)^{-1}$ as ${\cal F}_p^-\le {\cal E}_0$.

To prove \eqref{eq:A}  we introduce a new ${\cal G}$ valued 2-form 
$\B=r\A$ so that 
$$
\Ds_\L \B=0,\qquad \B|_{r=0}=\J.
$$
Considering the  components of the 2-form $\B$ it suffices to assume that $\B$ is a
${\cal G}$ valued function on ${\Bbb R}^{3+1}$, in fact on $\NN^-_{\tau_*}(p)$. 

Commuting\footnote{Recall that $\nab^{(\la)}$ is a gauge covariant derivative
acting on sections $S_r\to P\times_{Ad} {\cal G}$ 
and $\lap^{(\la)}$ is the corresponding  gauge Laplace-Beltrami operator 
on a standard sphere
$S_r$ of radius $r$. } the transport equation $\Ds_\L \B=0$ with 
$r^2\Delta^{(\la)}$ we obtain
$$
\Ds_\L (r^2\lap^{(\la)} \B)= r^2\, [\La_{\L}^{\,\, a}, \nab^{(\la)}_a \B] + 
r^2 \,\nab^{(\la)}_a 
[\La_{\L}^{\,\, a},
\B]
$$
We also have the equation 
$$
\Ds_\L (r \nab^{(\la)}_a \B)= r\, [\La_{\L a}, \B]
$$
We combine these equations into the system:
\begin{align}\label{eq:three-eqts}
&\Ds_\L \B=0,\\
&\Ds_\L (r \nab^{(\la)}_a \B)= r\, [\La_{\L a}, \B],\nn\\
&\Ds_\L (r^2\lap^{(\la)} \B)= 2 r^2\, [\La_{\L}^{\,\, a}, \nab^{(\la)}_a \B] + r^2 \,
[\nab_a^{(\la)}  \La_{\L}^{\,\,a}, \B].\nn
\end{align}
The first equation immediately implies that 
$$
\sup_{\NN^-_{\tau_*}(p)} |\B|\le |\J|,
$$
as the covariant derivative $\D_\a=\pa_\a+[\la.\c]$ is compatible with a scalar product on ${\cal G}$. 
We infer from the second equation that 
$$
\|\nab^{(\la)} \B\|_{L^2(S_r)}\le |\J|\, \sum_{a=1,2}
\int_0^r \left (\int_{S_\rho} |\La_{\L a}|^2 d\si_s\right)^{\f12} d\rho
$$
where $d\si_s$ is the are element of a $2$-dimensional sphere $S_\rho$ of
radius $\rho$. 

To treat the last equation in \eqref{eq:three-eqts}
we need to worry about the term $[\nab_a^{(\la)}  \La_{\L}^{\,\,a}, \B]$
which contains derivatives of $\La$. Recall that the flux only allows
us to estimate the tangential components of $\La$ and none 
if its derivatives. We get around this difficulty by expressing the Yang -Mills
equations $\D^\a \La_{\a\b}=0$ relative to the null frame $\L,\LLb, e_1, e_2$.
 In particular,
$\Ds^a \La_{\L a} + \Ds^\LLb \La_{\L \LLb}=0$. This in turn implies that 
$$
\nab_a^{(\la)} \La_{\L}^{\,\, a} = \f12 \Ds_\L \La_{\L \LLb} + \frac 1r \La_{\L\LLb}.
$$
Thus
$$
\Ds_\L (r^2\lap^{(\la)} \B-\f12 r^2 [\La_{\L\LLb}, \B])= 2 r^2\, [\La_{\L}^{\,\, a},
\nab^{(\la)}_a \B] 
$$
Therefore,
$$
\|\lap^{(\la)} \B\|_{L^1(S_r)}\le |\J|\,\left ( \f12\| \La_{\L\LLb}\|_{L^1(S_r)}+ \left
(\sum_{a=1,2}
\int_0^r \left (\int_{S_\rho} |\La_{\L a}|^2 d\si_s\right)^{\f12} d\rho\right)^2 \right)
$$
Integrating with respect to $r$ we obtain
$$
\|\lap^{(\la)} \B\|_{L^1(\NN^-_{\tau_*}(p))} \les |\J| \left (
\tau_*^{\frac 32} \| \La_{\L\LLb}\|_{L^2(\NN^-_{\tau_*}(p))}+ 
\tau_* \sum_{a=1,2} \|\La_{\L a}\|^2_{L^2(\NN^-_{\tau_*}(p))}\right ).
$$
and the result follows.
\end{proof}

\section{Applications to General Relativity}
In this section we specialize our results to Einstein vacuum space-times $(\M,\g)$: 
\begin{equation}\label{eq:vacuum}
\R_{\a\b}-\frac 12 \R\,\g_{\a\b}=0
\end{equation}
Equations \eqref{eq:vacuum} combined with the  Bianchi
 identities imply that
 the Riemann curvature tensor $\R_{\a\b\mu\nu}$ 
of an Einstein vacuum metric $\g$ satisfies a covariant wave equation 
$$
\Box_\g \R_{\a\b\mu\nu} = \left (\R\star \R\right)_{\a\b\mu\nu},
$$
where the quadratic term $\R\star \R$ is obtained by a 
contraction \footnote{These contractions 
result in a special structure of the quadratic term, crucial to 
the analysis in \cite{KR5} where
we investigate a breakdown criterion in General Relativity. The structure of this term is somewhat
analogous to the corresponding term in the Yang-Mills theory, see previous section.}
 of the curvature 
tensor $\R_{\a\b\mu\nu}$ with itself.
\begin{theorem}
\label{thm:KS-E}
Let $p$ be a point to the future of a space-like hypersurface $\Si$ in an Einstein vacuum
space-time $(\M,\g)$. We assume that assumptions ${\bf A1}, {\bf A2}$ are verified at $p$. 
Let $\A$ be a 4-tensor verifying,
 $$
\D_\L \A+\f12 \A \trch =0,\quad  s\A(p)=\J_p \qquad{\text on}\quad u=0
$$
where $\J_p$ is a fixed 4-tensor at $p$. Then the curvature tensor $\R_{\a\b\mu\nu}$ of $\g$
can be represented
by the following formula at point $p$ with $\Omega=\JJ^-(p)\cap \JJ^+(\Si)$,
\bea
4\pi \g(\R(p), \J_p)&=&-\int_{\Om}\de(u)\, \g(\A,\R\star \R) \label{eq:main-vac}\\
&-&\frac 12\int_{\Om}\de(u)\, \R(\c,\c, \LLb\, , \L)\#\R\#\A
+\int_{\Om} \de(u)\,\g\left(\big(\lap \A
+\ze_a \D_a \A \big),\R\right) \nn\\ &+&\f12\int_{\Om}
\de(u)\,\mu\,\g(\A,\R) + \int_{\Si} \left (\g(\A \de(u), \D_\T \R)-\g(\D_\T \left (\A\de(u)\right), \R)\right) \nn
\eea
where $\#$ denotes a contraction operation between tensors. 
The last term represents the contribution of the initial data on $\Si$. 
\end{theorem}
Representation \eqref{eq:main-vac} opens the  possibility of proving a pointwise bound
on the  curvature tensor in terms of initial data on $\Si$ and, as in the Yang-Mills
case, the  flux of curvature along the null boundary $\NN^-(p)$ of the set $\Omega$.
However, as opposed to the Yang-Mills equations on Minkowski background, where the
curvature flux is bounded by the $L^2$-norm of the curvature of initial data, no such a
priori bounds are available for the Einstein vacuum  equations. This suggests the use
of the $L^2$ based curvature norms to deduce  a breakdown criteria in General
Relativity, i.e. to show that  the space-time can be continued as long as such norms
remain finite.  For the Yang-Mills problem  in Minkowski
 space the underlying reason for
having 
 an a- priori bounds on the flux of curvature
is  due to the  presence of the  Killing vectorfield 
$\pr_t=\frac {\pa}{\pa t}$. In the case  of the Yang-Mills equations
 on  a smooth  curved background, such as in \cite{CS},  the result
remains true even though $\pr_t$ is no longer Killing; it suffices that
its deformation tensor is bounded. We  call such a vectorfield  {\it approximately
Killing}.

These considerations suggest  the following question. Assume that the
space-time 
$(\M,\g)$ possesses an {\it approximately Killing}, unit,  vectorfield $\T$,
 orthogonal to a space-like
Cauchy hypersurface
$\Si$, with    deformation tensor  $\pi(X,Y) = \g(\D_X \T, Y)$. We also 
 assume that the 
slices $\Si_t$ obtained by following integral curves of $\T$ from $\Si_0$ 
have constant mean curvature.
Can  the space-time  be extended as long as $\pi$ remains finite in the uniform norm?

The finiteness of the deformation tensor $\pi$ allows one to  control,
 via energy estimates
based on the  Bel-Robinson tensor, both  the $L^2$ norms
of the curvature  $\R$ along $\Si_t$ and the flux of curvature
along the  null boundaries
$\NN^-(p)$. The key step in the remaining analysis is 
 to derive a pointwise curvature based on the 
 representation formula \eqref{eq:main-vac}. In  \cite{KR5} we give
an affirmative answer to the question raised above by
 showing that the size of the 
region of validity of the formula 
\eqref{eq:main-vac} depends only on the  assumed
 $L^\infty$-  bounds on  
$\pi$ and  reasonable assumptions on the  initial data on $\Si$. Such estimates follow
from our work in
\cite{KR4}. More precisely we show that for a space-time metric $\g$ in the form 
$$
\g=-n^2 dt^2 + g_{ij} dx^i dx^j,
$$
where $n$ is the lapse function of the $t$ foliation and the vectorfield $\T$ is orthogonal to
$\Si_t$, the following result holds true.
  
\begin{theorem}\label{thm:main2} Assume that $(\M,\g)$ is a globally hyperbolic Einstein
vacuum space-time with $\Si_0$ a Cauchy hypersurface. 
Let the lapse function $n$ and the deformation
tensor $\pi$ of $\T$ satisfy 
$$
N_0^{-1}\le n\le N_0,\qquad \|\pi\|_{L^\infty}\le {\cal K}_0
$$
Assume  also that  $\M$
contains a future, compact set
${\cal D}\subset\M$  such that for any point $q\in {\cal D}^c$ the
 radius of injectivity of $\NN^-(q)$
 is at least $\de_0>0$.

Let $\Si$ be one of the slices of the $t$ foliation.
 Then assumptions ${\bf A1}, {\bf A2}$ of this paper are satisfied 
 for all points $p$ at distance $\le \de_*$
 from $\Si$, where $\de_*$ depends only on the Cauchy  data on $\Si_0$, $N_0$, ${\cal
K}_0$ and 
 $\de_0$. In particular, the  representation  formula \eqref{eq:main-vac} holds for all
such points.
\end{theorem}
\section{Open questions}
All the results of this paper have been derived under assumption
 ${\bf A2}$ which requires
that,  for any  point $p$,  the boundary $\NN^-(p)$ of the  causal past
 $\JJ^-(p)$ remains smooth 
at least until it reaches the 
 space-like   hypersurface  $\Si$. It is only  under this
 assumption that we
can guarantee that the  Kirchoff-Sobolev parametrix of Theorem \ref{thm:KS2} gives a
faithful representation  of a solution of the wave equation. We have 
already discussed the  two obstructions to smoothness of
$\NN^-(p)$: conjugate points of the  congruence of past 
directed null geodesics from $p$
and intersection of two distinct past  directed null geodesics from $p$. The second
obstruction can be easily   demonstrated on  a space-time 
$\M={\Bbb R}\times {\Bbb T}^3$ or $\M={\Bbb R}\times \Pi_a$ based on a flat torus 
${\Bbb T}^3$ or a flat cylinder $\Pi_a$ of width $a$. In those cases, however, 
an (exact) parametrix can be easily constructed by lifting the problem to the covering space 
${\Bbb R}\times {\Bbb R}^{3}$. On the other hand, the examples above are very special.  
Conjugate points can not be  
{\it removed} by  a simple\footnote{Conjugate points can be {\it  desingularized}
 however by lifting to the cotangent space } lifting  and the  quantity $\trch$,
which  features prominently 
in  our representation formulas, diverges to
$-\infty$  at a conjugate point. However, the same  focusing phenomenon shrinks the
volume of a {\it  conjugate point region}  thus leaving open a possibility that,
 perhaps,
with some additional assumptions on the structure and strength  of conjugate points,
  the  integral quantities appearing in 
our  representation formulas  remain 
finite.  It may thus be that the  Kirchoff-Sobolev parametrix remains valid 
even beyond
the region of formation of conjugate points. A good place to start 
  investigating this issue 
would be  product manifolds
$\M={\Bbb R}\times M$ with metrics of the form
$\g=-dt^2+g_{ij} dx^i dx^j$. A particularly interesting class to consider
  are product manifolds with  $M$ collapsing in the sense of Cheeger-Gromov.

\end{document}